\documentclass[a4paper]{article}

\usepackage{amsfonts, amsmath, amsthm, empheq, amssymb, bm, color,physics,enumerate}

\allowdisplaybreaks[4]

\newtheorem{lem}{Lemma}[section]
\newtheorem{thm}{Theorem}[section]
\newtheorem{prop}{Proposition}[section]
\newtheorem{cor}{Corollary}[section]
\newtheorem{rem}{Remark}[section]

\def\R{\mathbb R}

\def\De{\Delta}
\def\al{\alpha}
\def\be{\beta}

\def\ep{\epsilon}

\def\e{\mathrm e}
\def\f{\frac}

\def\Ga{\Gamma}

\def\la{\lambda}
\def\vp{\varphi}
\def\na{\nabla}

\def\Om{\Omega}
\def\ov{\overline}
\def\pa{\partial}
\def\Re{\mathrm{\,Re\,}}
\def\Im{\mathrm{\,Im\,}}

\title{\Large\bf\boldmath
Carleman estimate for the Schr\"{o}dinger equation and application to magnetic inverse problems}

\author{\large Xinchi HUANG, Yavar KIAN, \'Eric SOCCORSI and Masahiro YAMAMOTO}

\begin{document}
\maketitle


\section{Introduction}
\label{sec-multi}
Let $T \in (0,+\infty)$ and let $\Om$ be a bounded domain 
in $\R^n$, $n\in \Bbb{N}:=\{1,2,\ldots \}$, with smooth boundary $\pa\Om$. We consider the following initial-boundary value problem for the magnetic Schr\"{o}dinger equation
\begin{equation}
\label{sys:ori}
\left\{
\begin{aligned}
& -i\pa_t u - \De_{A_0} u + \rho_0 u = 0 \qquad && \text{in}\ \Om\times (0,T):=Q, \\
& u = g \qquad && \text{on}\ \pa\Om\times (0,T):=\Sigma, \\
& u(\cdot,0) = u_0&& \text{in}\ \Om, 
\end{aligned}
\right.
\end{equation}
with initial state $u_0$ and non-homogeneous Dirichlet boundary condition $g$. Here $\rho_0 : \Om \to \Bbb{C}$ is a complex-valued electric potential and 
\begin{equation}
\label{def-deltaA}
\De_{A_0} := (\na + iA_0)\!\cdot\! (\na + iA_0) = \De + 2iA_0\!\cdot\!\na + i(\na\!\cdot\! A_0) - \abs{A_0}^2
\end{equation}
denotes the magnetic Laplace operator associated with the magnetic vector potential
$A_0 : \Omega \to \R^n$. In the particular case where $n=3$, the magnetic field induced by the magnetic potential vector $A_0$ reads ${\rm curl}\ A_0 :=\nabla \cross A_0$. 

In the present paper we examine stability issue in the inverse problem of determining the electromagnetic potential $(\rho_0,A_0)$ from a finite number of partial Neumann boundary measurements over the entire time-span of the solution to \eqref{sys:ori}, by $n+1$ times suitably changing the initial state $u_0$.

There are numerous papers available in the mathematical literature dealing with inverse coefficients problems from knowledge of the Dirichlet-to-Neumann (DN) map. But in the particular case of the magnetic Schr\"odinger equation, the DN map $\Lambda_{A_0}$ is invariant under gauge transformation of $A_0$, i.e. $\Lambda_{A_0 + \nabla \psi}= \Lambda_{A_0}$ for all $\psi \in C^1(\overline{\Omega})$ such that $\psi_{| \partial \Om}=0$, see e.g. \cite{Eskin1}.
Therefore the magnetic potential vector cannot be uniquely determined by the DN map and the best we can expect from knowledge of $\Lambda_{A_0}$ is uniqueness modulo gauge transform of $A_0$. However it is well known that the magnetic field $d A_0$, i.e. the exterior derivative of $A_0$ defined as the 1-form $\sum_{j=1}^n (A_0)_j d x_j$ (if $n=3$ then $d A_0$ is generated by the curl of $A_0$) is invariant under gauge transformation of $A_0$. As a matter of fact it was proved in \cite{Sun} that the DN map uniquely determines the magnetic field provided the underlying magnetic vector potential is sufficiently small in a suitable class. The smallness assumption was removed in \cite{NSU} for $C^{\infty}$ magnetic vector potentials and later on this smoothness assumption was weakened to $C^1$ in \cite{Tol} and to Dini continuous in \cite{Sal}. 
In \cite{Eskin2}, the author proves that the electromagnetic potential of the Schr\"odinger equation in domains with several obstacles, is uniquely defined by the DN map.
In \cite{BeCh10} the magnetic field is stably retrieved by the dynamical DN map. The uniqueness and stability issues for 
time-dependent electromagnetic potentials of the Schr\"odinger equation are addressed in \cite{Eskin3} and \cite{KS18}, respectively.  All the above cited results were obtained with the full DN map, which is made of measurements of the solution taken on the whole boundary. The uniqueness problem by a local DN map was solved in \cite{FKSU} and it was shown in \cite{Tzou} that the magnetic field depends stably on the DN map measured on any sub-boundary that is slightly larger than half the boundary. This result was extended in \cite{BenJ} to arbitrary small sub-boundaries provided the magnetic potential is known in the vicinity of the boundary.

Notice that infinitely many boundary observations of the solution to the magnetic Schr\"odinger equation were needed in all the above mentioned articles in order to define the DN map. By contrast, the time independent and real valued electric potential in a Schr\"odinger equation was stably retrieved from a single boundary measurement in \cite{BP02, MAR}. In these two papers, the observation zone fulfills a geometric condition related to geometric optics condition insuring observability. This geometric condition was relaxed in \cite{BC} upon assuming that the electrostatic potential is known near the boundary. In \cite{CS}, the space varying part of the divergence free $n$ dimensional magnetic potential was reconstructed by $n$ partial Neumann data, by changing the initial state of the Schr\"odinger equation $n$ times suitably.

In the present article we aim for stable determination of the electromagnetic potential $(\rho_0,A_0)$ in \eqref{sys:ori} through $n+1$ partial Neumann observations, by means of a Carleman estimate. We refer to \cite{BP02,Albano,Tataru,YY10} for actual examples of Carleman inequalities for the Schr\"odinger equation. The idea of using Carleman estimates for solving inverse coefficient problems was first introduced by Bugkheim and Klibanov in \cite{BK:81}. Since then, this technique has then been successfully applied by numerous authors to various types (parabolic, hyperbolic, elasticity, Maxwell, etc.) inverse coefficients problems in bounded domains, see e.g. \cite{ImaYama01, K:92, KT, BCS, Y99} and references therein (and more recently it was adapted to the reconstruction of non compactly supported unknown coefficients in \cite{KPS1,KPS2, BKS, CLS}). More specifically, in the framework of the Schr\"odinger equation, the authors of \cite{BP02, CS, YY10} use a Carleman inequality on the extended domain 
$\Omega \times (-T,T)$ in order to avoid observation data at $t=0$ over $\Om$, appearing in Carleman estimates on $Q$. This imposes that the solution $u$ to \eqref{sys:ori}, extended to $\Om \times (-T,T)$ by setting either $u(x,t)=\overline{u(x,-t)}$ or $u(x,t)=-\overline{u(x,-t)}$ for a.e. $(x,t) \in \Om \times (-T,0)$, depending on whether the initial state $u_0$ is taken real-valued or purely imaginary, be a solution to the Schr\"odinger equation in $\Om \times (-T,T)$. It follows readily from the above time-symmetrization $u(\cdot,t)=\pm \overline{u(\cdot,-t)}$ that
\begin{equation}
\label{sy-1}
(-i \partial_t - \Delta_{-A_0} + \ov{\rho_0}) u(\cdot,t) = \pm (\ov{(-i \partial_t - \Delta_{A_0} + \rho_0) u})(\cdot,-t),\
t \in (-T,0),
\end{equation}
and hence that $u$ is solution to the Schr\"odinger equation in $\Om \times (-T,T)$ iff $(-A_0,\ov{\rho_0})=(A_0,\rho_0)$, i.e. $A_0=0$ and $\rho_0 \in \R$ (this is precisely the situation examined in \cite{BP02}, where Lipschitz stable reconstruction of the real-valued electrostatic potential $\rho_0$ is derived in absence of a magnetic potential), in which case the right hand side of \eqref{sy-1} is zero. As a conclusion, the time-symmetrization method implemented in \cite{BP02} does not work in presence of a non-zero time-independent magnetic potential vector $A_0$ (notice that this is no longer true for odd time-dependent magnetic potentials: Indeed, when $\rho_0$ and $A_0$ depend on $(x,t)$ then \eqref{sy-1} reads
$(-i \partial_t - \Delta_{-A_0(\cdot,t)} + \ov{\rho_0(\cdot,t)}) u(\cdot,t) = \pm (\ov{(-i \partial_t - \Delta_{A_0(\cdot,-t)} + \rho_0(\cdot,-t)) u})(\cdot,-t)$ for a.e. $t \in (-T,0)$, 
so the extended solution $u$ fulfills the magnetic Schr\"odinger equation in $\Om \times (-T,T)$ iff we have
$(-A_0(\cdot,t),\ov{\rho_0(\cdot,t)})=(A_0(\cdot,-t),\rho_0(\cdot,-t))$, which corresponds to the framework of \cite{CS,YY10}). Therefore, in contrast with \cite{BP02,CS,YY10}, we cannot symmetrize the solution to \eqref{sys:ori} with respect to the time-variable in the framework of in this paper. As a consequence we need a modified global Carleman estimate for the Schr\"odinger operator in $Q$, as compared to the ones of \cite{BP02,YY10} that are established in $\Om \times (-T,T)$, in order to adapt the Bukhgeim-Klibanov method to the "stationary magnetic" Schr\"odinger equation under investigation here. We shall actually prove the following three stability results for the inverse problem under inverstigation.
\begin{enumerate}[i)]
\item {\it Case 1}: Assuming that $A_0$ is known, we stably determine the complex-valued electric potential $\rho_0$ from a single partial boundary measurement over the entire time span of the normal derivative
of the solution $u$ to \eqref{sys:ori}, measured on a sub-boundary $\Ga_0\subset \pa\Om$. The result is valid for any two electrostatic potentials with difference $\rho$, whose imaginary part of the logarithmic gradient $\nabla \ln (\rho^{-1} \ov{\rho} )$ is uniformly bounded in $\Om$, see condition \eqref{c-el} below.
\item {\it Case 2}: We prove simultaneous stable reconstruction of the magnetic vector potential $A_0$ (together with its divergence $\nabla \cdot A_0$) and the complex-valued electric potential $\rho_0$, through $n+1$ partial Neumann observations of the solution, obtained by changing $n+1$ times the initial condition $u_0$ suitably. This is provided the logarithmic gradient of the difference of the electromagnetic potentials is uniformly bounded in $\Om$, see assumptions \eqref{h-em-a}, \eqref{h-em-b} and \eqref{h-em-c}.
\item {\it Case 3}: Assuming that $\rho_0$ and the strength $\abs{A_0}$ of the magnetic potential vector are known, we stably retrieve the direction of $A_0$ (together with the divergence), from $n+1$ partial Neumann data. In contrast with the two above results, there is no additional condition of the type of \eqref{c-el} or \eqref{h-em-a}-\eqref{h-em-c}, imposed on the magnetic vector potential for this result to hold.
\end{enumerate}

Our first claim (see Theorem \ref{thm:stb} below) extends the stability results of \cite{BP02} to the case of complex-valued electrostatic potentials. We refer to  \cite{KE}[Part 2, Section 14, Appendix B] for the physical relevance of complex-valued electric potentials appearing in the Schr\"odinger equation. Moreover, to the best of our knowledge, the second and third claims (see Theorems \ref{thm:stb2} and \ref{thm:stb3}) are the only stability results by finitely many local Neumann data, for stationary magnetic potential vectors of the Schr\"odinger equation, available in the mathematical literature.

\subsection{Notations}
Throughout this text $x:=(x_1,\ldots,x_n)$ denotes a generic point of $\Om \subset \R^n$ and we use the following notations: $\pa_i := \frac{\pa}{\pa x_i}$ for $i=1,\ldots,n$. We put $\pa^2_{i,j}:=\pa_i \pa_j$ for $i,j=1,\ldots,n$ and as usual we write $\pa_i^2$ instead of $\pa_{ii}^2$. Similarly, we write $\pa_t := \frac{\pa}{\pa t}$ and $\pa_\nu u = \f{\pa u}{\pa \nu} := \na u \cdot \nu$, where $\nu$ denotes the outward normal vector to the boundary $\pa\Om$ and $\nabla$ is the gradient operator with respect to the space variable $x$. The symbol $\cdot$ denotes the Euclidian scalar product in $\R^n$ and $\nabla \cdot$ stands for the divergence operator.

Let us now introduce the following functional spaces. For $X$, a manifold, we set
$$ H^{r,s}(X \times (0,T)):=L^2(0,T;H^r(X)) \cap H^s(0,T;L^2(X)).$$ 
When $X=\Omega$, we write 
$H^{r,s}(Q)=L^2(0,T;H^r(\Omega)) \cap H^s(0,T;L^2(\Omega))$ instead of $H^{r,s}(\Om \times (0,T))$ and for $X=\partial \Omega$, we write
$H^{r,s}(\Sigma)=L^2(0,T;H^r(\partial \Omega)) \cap H^s(0,T;L^2(\partial \Omega))$ instead of $H^{r,s}(\pa \Om \times (0,T))$.

\subsection{Existence and uniqueness results}
\label{sec-eu}

Our first result is as follows.

\begin{prop}
\label{pr-eu}
For $M \in (0,+\infty)$, let $A_0 \in W^{1,\infty}(\Omega,\R^n)$ and $\rho_0 \in W^{2,\infty}(\Omega,\Bbb{C})$ satisfy
\begin{equation}
\label{sc}
\norm{A_0}_{W^{1,\infty}(\Omega)} +  \norm{\rho_0}_{W^{2,\infty}(\Omega)} \leq M.
\end{equation}
Then, for all $g \in H^{\frac{7}{2},\frac{7}{4}}(\Sigma)$ and all $w_0 \in H^2(\Omega)$ obeying
\begin{equation}
\label{cc}
g(\cdot,0)=u_0\ \mbox{on}\ \partial \Omega,
\end{equation}
there exists a unique solution $u \in H^{2,1}(Q)$ to \eqref{sys:ori}. Moreover we have the estimate
\begin{equation}
\label{es}
\norm{u}_{H^{2,1}(Q)}  \leq   C \left( \norm{u_0}_{H^2(\Omega)} + \norm{g}_{H^{\frac{7}{2},\frac{7}{4}}(\Sigma)} \right),
\end{equation}
where $C$ is a positive constant depending only on $T$, $\Omega$ and $M$. 
\end{prop}

As a corollary we have the following improved regularity result.

\begin{cor}
\label{c-eu}
For $M \in (0,+\infty)$, let $A_0 \in W^{3,\infty}(\Omega,\R^n)$ and $\rho_0 \in W^{2,\infty}(\Omega,\Bbb{C})$ fulfill
\begin{equation}
\label{c-sc}
\norm{A_0}_{W^{3,\infty}(\Omega)} +  \norm{\rho_0}_{W^{2,\infty}(\Omega)} \leq M.
\end{equation}
Pick $g \in H^{\frac{11}{2},\frac{11}{4}}(\Sigma)$ and $u_0 \in H^4(\Omega)$ such that
\begin{equation}
\label{c-cc}
\partial_t^k g(\cdot,0) = \left(- i (-\Delta_{A_0}+\rho_0) \right)^k u_0\ \mbox{on}\ \partial \Omega\ \mbox{for}\ k=0,1.
\end{equation}
Then there exists a unique solution $u \in H^{2}(0,T;L^2(\Omega)) \cap H^1(0,T;H^2(\Omega))$ to \eqref{sys:ori}, satisfying
\begin{equation}
\label{c-es}
\norm{u}_{H^{2,1}(Q)}  + \norm{\partial_t u}_{H^{2,1}(Q)}\leq   C \left( \norm{u_0}_{H^4(\Omega)} + \norm{g}_{H^{\frac{11}{2},\frac{11}{4}}(\Sigma)} \right),
\end{equation}
for some positive constant $C$ depending only on $T$, $\Omega$ and $M$. 
\end{cor}

\subsection{Inverse problem: main results}

Given $M \in (0,+\infty)$ and $q_0 \in W^{2,\infty}(\Omega,\R)$ we define the set of admissible unknown electric potentials as
$$ \mathcal{Q}_M(q_0) : = \{ q \in W^{2,\infty}(\Omega,\Bbb{C}):\ \norm{q}_{W^{2,\infty}(\Omega)} \leq M\ \mbox{and}\ q = q_0\ \mbox{on}\ \partial \Omega \}. $$

We first address the inverse problem of recovering the complex-valued electrostatic potential when the magnetic vector potential is known.

\begin{thm}
\label{thm:stb}
Let $g$ and $u_0$ be as in Corollary \ref{c-eu}. Moreover, we suppose that
\begin{equation}
\label{ass:u_0}
\exists r_0 \in (0,+\infty),\ \abs{u_0(x)} \ge r_0,\qquad x\in\Om.
\end{equation}
For $M \in (0,+\infty)$ and $\rho_0 \in W^{2,\infty}(\Omega,\Bbb{C})$, let $\rho_j \in \mathcal{Q}_M(\rho_0)$, $j=1,2$, fulfill the following condition a.e. in $\Omega$:
\begin{equation}
\label{c-el}
\abs{\Im \left( \overline{(\rho_1-\rho_2)} \nabla (\rho_1-\rho_2 ) \right)}\le M \abs{\rho_1-\rho_2}^2.
\end{equation}
Then there exist a nonempty sub-boundary $\Ga_0\subset\pa\Om$ and a positive constant $C$ depending only on $\Omega$, $T$, $M$ and $\rho_0$, such that
\begin{equation*}
\norm{\rho_1-\rho_2}_{L^2(\Om)}\le C \norm{\pa_\nu \pa_t (u_1-u_2)}_{L^2(\Ga_0\times (0,T))}.
\end{equation*}
Here $u_j$, for $j=1,2$, is the $H^2(0,T;L^2(\Om)) \cap H^1(0,T;H^2(\Om))$-solution to \eqref{sys:ori} associated with electric potential $\rho_j$ and uniformly zero magnetic potential, which is given by Corollary \ref{c-eu}. 
\end{thm}

Let us briefly comment on Theorem \ref{thm:stb}:
\begin{enumerate}[a)]
\item The assumption \eqref{ass:u_0} allows for a far more flexible choice of initial input $u_0$ than in \cite{BP02, CS}, where it is required to be either real-valued or purely imaginary.
\item The condition \eqref{c-el} holds true provided either of the real or imaginary parts of the electrostatic potential, is known.  
Therefore Theorem \ref{thm:stb} extends the stability result of \cite{BP02}. 
\item Arguing as in \cite{BP02}, we can prove at the expense of higher regularity on the coefficients and data of the magnetic Schr\"odinger equation, that the following double-sided stability inequality
$$
\norm{\rho_1-\rho_2}_{H_0^1(\Om)}\le C_1\norm{\pa_\nu \pa_t (u_1-u_2)}_{L^2(\Ga_0\times (0,T))}\le C_2\norm{\rho_1-\rho_2}_{H_0^1(\Om)},
$$
holds for two positive constants $C_1$ and $C_2$.
\item There are actual classes of complex-valued electrostatic potentials fulfilling condition \eqref{c-el}. For instance, this is the case of $\mathcal{E}_a := \{ \rho(x) = a + \delta \langle x \rangle,\ \delta \in \Bbb{C} \}$, where $a \in \Bbb{C}$ is arbitrary and $\langle x \rangle := (1+ \abs{x}^2)^{\f{1}{2}}$ for $x \in \R^n$. Indeed, for any $\rho_j(x) = a + \delta_j \langle x \rangle \in \mathcal{E}_a$, $j=1,2$, it holds true that 
$$\abs{\nabla (\rho_1-\rho_2)(x)}= \abs{\delta_1-\delta_2} \f{\abs{x}}{ \langle x \rangle} \le \abs{\delta_1-\delta_2} \le\abs{(\rho_1-\rho_2)(x)},\ x \in \Om. $$  

\end{enumerate}

Next we consider the inverse problem of determining the electromagnetic potential $(A_0,\rho_0)$.  For $A_0 \in W^{3,\infty}(\Omega,\R^n)$, we define the set of admissible unknown magnetic potentials as
$$ 
\mathcal{A}_M(A_0) := \{ A \in W^{3,\infty}(\Omega,\R^n):\  \norm{A}_{W^{3,\infty}(\Omega)} \le M\ \mbox{and}\ A=A_0,\  \nabla \cdot A = \nabla \cdot A_0 \ \mbox{on}\ \partial \Omega \}. 
$$

\begin{thm}
\label{thm:stb2}
Pick $M \in (0,+\infty)$, $\rho_0 \in W^{2,\infty}(\Omega,\Bbb{C})$ and $A_0 \in W^{3,\infty}(\Omega,\R^n)$, and for $j=1,2$, let $\rho_j \in \mathcal{Q}_M(\rho_0)$ and $A_j \in \mathcal{A}_M(A_0)$ fulfill the three following conditions a.e. in $\Omega$:
\begin{eqnarray}
& & \abs{\nabla (\rho_1-\rho_2)} \le M \abs{(\rho_1-\rho_2)}, \label{h-em-a} \\
& & \max_{i=1,\ldots,n} \sum_{j=1}^n \abs{\partial_{x_i} (A_1-A_2)_j}  \le  M \abs{(A_1-A_2)}, \label{h-em-b} \\
& & \abs{\nabla ( \nabla \cdot (A_1-A_2) )} \le M \abs{\nabla \cdot (A_1-A_2)}. \label{h-em-c}
\end{eqnarray}
Then for all $g \in H^{\frac{11}{2},\frac{11}{4}}(\Sigma)$, there exists a set of initial states $\{u_0^k\}_{k=0}^n\subset (W^{1,\infty}(\Om, \Bbb{C}) \cap H^4(\Om,\Bbb{C}))^{n+1}$ obeying the condition \eqref{c-cc}, such that the stability inequality
$$
\norm{\rho_1-\rho_2}_{L^2(\Om)} + \norm{A_1-A_2}_{L^2(\Om)} + \norm{\na\!\cdot\! A_1 - \na\!\cdot\! A_2}_{L^2(\Om)} \le C\sum_{k=0}^n \norm{\pa_\nu \pa_t (u_1^k-u_2^k)}_{L^2(\Ga_0\times (0,T))}
$$
holds for some nonempty sub-boundary $\Gamma_0\subset\pa\Om$ and some positive constant $C$, depending only on $T$, $\Om$ and $M$. Here $u_j^k$, $j=1,2$ and $k=0,\ldots,n$, denotes the $H^{2}(0,T;L^2(\Omega)) \cap H^1(0,T;H^2(\Omega))$-solution to \eqref{sys:ori} with initial state $u_0^k$ and electromagnetic potential $(\rho_j,A_j)$.
\end{thm}

Actual examples of classes of electromagnetic potentials fulfilling conditions \eqref{h-em-a}-\eqref{h-em-c} can be built in the same fashion as in Point d) following Theorem \ref{thm:stb}. 

Finally, we consider the inverse problem of determining the direction of the magnetic potential vector when its strength, together with the electric potential, are known.

\begin{thm}
\label{thm:stb3}
For $M \in (0,+\infty)$ and $A_0 \in W^{3,\infty}(\Omega,\R^n)$, let $A_j \in \mathcal{A}_M(A_0)$, $j=1,2$, be such that
\begin{equation}
\label{h-em2}
\abs{A_1(x)} = \abs{A_2(x)},\ x \in \Om.
\end{equation}
Then for all $g \in H^{\frac{11}{2},\frac{11}{4}}(\Sigma)$, there exists a set of initial states $\{u_0^k\}_{k=0}^n\subset (W^{1,\infty}(\Om, \Bbb{C}) \cap H^4(\Om,\Bbb{C}))^{n+1}$ fulfilling \eqref{c-cc}, such that we have 
$$
\norm{A_1-A_2}_{L^2(\Om)}  + \norm{\nabla \cdot A_1- \nabla \cdot A_2}_{L^2(\Om)} \le C\sum_{k=0}^n \norm{\pa_\nu \pa_t (u_1^k-u_2^k)}_{L^2(\Ga_0\times (0,T))}
$$
for some nonempty sub-boundary $\Gamma_0\subset\pa\Om$ and some positive constant $C$, depending only on $T$, $\Om$ and $M$. Here $u_j^k$, for $j=1,2$ and $k=0,\ldots,n$, is the $H^{2}(0,T;L^2(\Omega)) \cap H^1(0,T;H^2(\Omega))$-solution to \eqref{sys:ori} with initial state $u_0^k$, electric potential $\rho_0=0$ and magnetic potential $A_0=A_j$.
\end{thm}

If the divergence of the magnetic vector potentials is known, in such a way that $\nabla \cdot (A_1-A_2)=0$ everywhere in $\Om$, then it is easy to see from the derivation of Theorem \ref{thm:stb3}, given in Subsection \ref{sec-stb3}, that the above stability inequality remains valid with only  $n$ local boundary measurements. We point that such a result is optimal in the sense that we recover the $n$ components of an unknown vector-valued function by exactly $n$ local boundary measurements.

\subsection{Overview}
The paper is organized as follows. In Section \ref{sec-pr-eu} we study the forward problem associated with \eqref{sys:ori}: We prove Proposition \ref{pr-eu} and Corollary \ref{c-eu}.
Section \ref{sec-car} is devoted to the derivation of a Carleman estimate for the Schr\"odinger equation in $Q$, needed by the analysis of the inverse problems under examiniation. Finally, Section \ref{sec-pr} contains the proof of Theorems \ref{thm:stb}, \ref{thm:stb2} and \ref{thm:stb3}.


\section{Analysis of the direct problem}
\label{sec-pr-eu}
For $A_0 \in W^{1,\infty}(\Om,\R)$ we denote by $-\Delta_{A_0}^D$ the selfadjoint operator generated in $L^2(\Omega)$ by the closed symmetric form
$$ a(u,v) := \int_{\Omega} (\nabla + i A_0) u(x) \cdot \overline{(\nabla + i A_0) v(x)} dx,\ u,v \in H_0^1(\Omega). $$
It is well known that the Dirichlet Laplacian $-\Delta_{A_0}^D$ acts on his domain $H_0^1(\Omega) \cap H^2(\Omega)$ as the operator $-\Delta_{A_0}$. 

Let $A_0$ and $\rho_0$ be the same as in Proposition \ref{pr-eu}.
Then, upon applying \cite{CKS}[Lemma 2.1] (with $X=L^2(\Omega)$, $U=i \Delta_{A_0}^D$ and $B(t)=- i \rho_0$ for any $t \in [0,T]$), we obtain the :

\begin{lem}
\label{lm-cp}
Let $v_0 \in H_0^1(\Omega) \cap H^2(\Omega)$ and let $f \in H^{0,1}(Q)$. Then there is a unique solution
$v \in C([0,T],H_0^1(\Omega) \cap H^2(\Omega))\cap C^1([0,T],L^2(\Omega))$ to 
\begin{equation}
\label{cp1}
\left\{
\begin{array}{ll}
(-i \partial_t  - \Delta_{A_0}^D + \rho_0 ) v  = f \\
v(\cdot,0)=v_0, 
\end{array}
\right.
\end{equation}
fulfilling
\begin{equation}
\label{cp2}
\norm{v(\cdot,t)}_{H^2(\Omega)}+ \norm{\partial_t v(\cdot,t)}_{L^2(\Omega))} \leq C \left( \norm{v_0}_{H^2(\Omega)}+ \norm{f}_{H^{0,1}(Q)} \right),
\end{equation}
uniformly in $t \in [0,T]$. Here $C$ is some positive constant depending only on $\Omega$, $T$ and $M$.
\end{lem}

\subsection{Proof of Proposition \ref{pr-eu}}
Since $g \in H^{\frac{7}{2},\frac{7}{4}}(\Sigma)$ then \cite{LM2}[Section 4, Theorem 2.3] yields existence of $G \in H^{4,2}(Q)$ such that $G=g$ on $\Sigma$ and 
\begin{equation}
\label{z0}
\norm{G}_{H^{4,2}(Q)} \leq C \norm{g}_{H^{\frac{7}{2},\frac{7}{4}}(\Sigma)},
\end{equation}
for some positive constant $C$ depending only on $T$ and $\Omega$. Therefore $w$ solves \eqref{sys:ori}
iff the function $\tilde{w}:=w-G$ is solution to
\begin{equation}
\label{st}
\left\{ \begin{array}{ll}
(-i \partial_t - \Delta_{A_0} + \rho_0) \tilde{w} = f & \mbox{in}\ Q \\
\tilde{w} = 0 & \mbox{on}\ \Sigma \\
\tilde{w}(\cdot,0) = \tilde{w}_0 & \mbox{in}\ \Omega, \end{array} \right.
\end{equation}
with $f:=-(-i \partial_t - \Delta_{A_0} + \rho_0) G$ and  $ \tilde{w}_0:=w_0 -G(\cdot,0)$. Next, as $G \in H^{4,2}(Q)$ yields $\partial_t G \in H^{2,1}(Q)$ in virtue of , then it is apparent that $f \in H^{0,1}(Q)$ and
\begin{equation}
\label{z1}
\norm{f}_{H^{0,1}(Q)} \leq C \left( \norm{G}_{H^{2,1}(Q)} + \norm{\partial_t G}_{H^{2,1}(Q)} \right) \leq C  \norm{G}_{H^{4,2}(Q)}.
\end{equation}
Further we have $G(\cdot,0) \in H^{3}(\Omega)$ by \cite{LM2}[Section 4, Theorem 2.1], with 
$ \norm{G(\cdot,0)}_{H^3(\Omega)} \leq C \norm{G}_{H^{4,2}(Q)}$,
whence $\tilde{w}_0 \in H^2(\Omega)$ and
\begin{equation}
\label{z2}
\norm{\tilde{w}_0}_{H^2(\Omega)} \leq C \left( \norm{w_0}_{H^2(\Omega)} + \norm{G}_{H^{4,2}(Q)} \right).
\end{equation}
Moreover, since $\tilde{w}_0=0$ on $\partial \Omega$, by \eqref{cc}, then we have $\tilde{w}_0 \in H_0^1(\Omega) \cap H^2(\Omega)$ so we may apply Lemma \ref{lm-cp} to \eqref{st}. We get that there is a unique solution $\tilde{w} \in H^{2,1}(Q)$ to
\eqref{st}, such that
$$
\norm{\tilde{w}}_{H^{2,1}(Q)} \leq C \left( \norm{\tilde{w}_0}_{H^2(\Omega)} + \norm{f}_{H^{0,1}(Q)} \right).
$$
Finally, putting this together with the estimate $\norm{w}_{H^{2,1}(Q)} \leq  \left( \norm{\tilde{w}}_{H^{2,1}(Q)} + \norm{G}_{H^{2,1}(Q)} \right)$, \eqref{z0} and \eqref{z1}, we obtain \eqref{es}.

\subsection{Proof of Corollary \ref{c-eu}}
Put $z:=\partial_t w$, where $w$ is the $H^{2,1}(Q)$-solution to \eqref{sys:ori}, given by Proposition \ref{pr-eu}. Then we have
\begin{equation}
\label{c-s}
\left\{ \begin{array}{ll}
(-i \partial_t - \Delta_{A_0} + \rho_0) z = 0 & \mbox{in}\ Q\\
z = \partial_t g & \mbox{on}\ \Sigma \\
z(\cdot,0) = z_0 & \mbox{in}\ \Omega, \end{array} \right.
\end{equation}
where $z_0 := -i (-\Delta_{A_0}+\rho_0) w_0 \in H^2(\Omega)$. 
Moreover we have $\partial_t g \in H^{\frac{7}{2},\frac{7}{4}}(\Sigma)$ from \cite{LM2}[Section 4, Proposition 2.3], with
$$
\norm{\partial_t g}_{H^{\frac{7}{2},\frac{7}{4}}(\Sigma)} \leq C \norm{g}_{H^{\frac{11}{2},\frac{11}{4}}(\Sigma)}.
$$
Therefore, in light of \eqref{c-s} and the compatibility condition \eqref{c-cc} with $k=1$, we infer from Proposition \ref{pr-eu} that $z \in H^{2,1}(Q)$ satisfies
\begin{eqnarray*}
\norm{z}_{H^{2,1}(Q)} & \leq & C \left(  \norm{z_0}_{H^2(\Omega)} + \norm{\partial_t g}_{H^{\frac{7}{2},\frac{7}{4}}(\Sigma)} \right) \nonumber \\
& \leq & C \left(  \norm{w_0}_{H^4(\Omega)} +  \norm{g}_{H^{\frac{11}{2},\frac{11}{4}}(\Sigma)} \right). 
\end{eqnarray*}
The desired result follows readily from this and \eqref{es}.

\section{Global Carleman estimate}
\label{sec-car}
In this section, we establish a global Carleman estimate for the main part of Schr\"{o}dinger operator
\begin{equation}
\label{equ-gov}
L:=-i \pa_t  - \De  
\end{equation}
acting in $Q=\Om \times (0,T)$. 
Carleman estimates for the Schr\"{o}dinger operator in domains centered around $t=0$ such as $\Om \times (-T,T)$ were derived in \cite{YY10} with a regular weight function and in \cite{BP02} with a symmetric singular weight function. However, since the solution $u$ to \eqref{sys:ori} cannot be time-symmetrized in the framework of this paper, we need to establish a Carleman estimate for the operator $L$ in $Q$.

To this end, we assume in the entire section that $u \in L^2(0,T; H^1_0(\Om))$ and $Lu\in L^2(Q)$. Notice for further use that $\f{\pa u}{\pa \nu}\in L^2(\Sigma)$. Next we put $w:=ue^{s\al}$, where $s \in (0,+\infty)$ and $\al$ is a real-valued smooth function we shall make precise further, and set
\begin{align*}
Rw:=e^{s\al}Lu = e^{s\al}L(we^{-s\al}) = is(\pa_t \al)w + R_1 w + R_2 w = is(\pa_t \al)w + R_3 w,
\end{align*}
with
\begin{align}
\nonumber
&R_1 w := -i\pa_t w - \De w - s^2|\na\al|^2 w,\\
\label{def:R}
&R_2 w := 2s\na\al\!\cdot\!\na w + s(\De\al) w,\\
\nonumber
&R_3 w := R w - is(\pa_t \al)w = R_1 w + R_2 w.
\end{align}
Since the function $w$ is complex-valued, we denote by $w_{{\rm re}}$ its real part and by $w_{{\rm im}}$ its imaginary part, in such a way that $w = w_{{\rm re}} + iw_{{\rm im}}$.
Similarly we decompose each $R_j w$, for $j=1,2,3$, into the sum
$$
R_j w =P_j w + iQ_j w,
$$
where
\begin{align*}
P_1 w := \pa_t w_{{\rm im}} - \De w_{{\rm re}} -s^2|\na\al|^2 w_{{\rm re}},\ & Q_1 w := - \pa_t w_{{\rm re}} - \De w_{{\rm im}} -s^2|\na\al|^2 w_{{\rm im}}, \\
P_2 w := 2s\na\al\!\cdot\!\na w_{{\rm re}} + s (\De\al) w_{{\rm re}},\ & Q_2 w := 2s\na\al\!\cdot\!\na w_{{\rm im}} + s (\De\al) w_{{\rm im}}, \\
P_3 w := \mathrm{Re}(R w) + s(\pa_t\al) w_{{\rm im}},\ & Q_3 w := \mathrm{Im}(R w) - s(\pa_t\al) w_{{\rm re}}.\end{align*}
As we are aiming for computing $\abs{R_3 w}^2$ and since
\begin{align}
\label{eq1}
\abs{R_3 w}^2 = \sum_{j=1}^2 \abs{R_j w}^2  + 2 \Re \left( (R_1 w) \ov{R_2 w} \right) =  \sum_{j=1}^2 \abs{R_j w}^2
+ 2 (P_1 w) P_2 w + 2 (Q_1 w) Q_2 w,
\end{align}
we start by expanding the two last terms in the right hand side of \eqref{eq1}. We get that
\begin{align*}
2(P_1 w, P_2 w)_{L^2(Q)} &= \int_Q 4s (\na\al\!\cdot\!\na w_{{\rm re}}) (\pa_t w_{{\rm im}}) dxdt + \int_{Q} 2s(\De\al) w_{{\rm re}}(\pa_t w_{{\rm im}}) dxdt - \int_{Q} 4s(\na\al\!\cdot\!\na w_{{\rm re}})\De w_{{\rm re}} dxdt \\
&- \int_{Q} 2s (\De\al) w_{{\rm re}} \De w_{{\rm re}} dxdt - \int_{Q} 4s^3 |\na\al|^2(\na\al\!\cdot\!\na w_{{\rm re}})w_{{\rm re}} dxdt - \int_{Q} 2s^3 |\na\al|^2(\De\al) |w_{{\rm re}}|^2 dxdt \\
&=: \sum_{k=1}^6 I_k
\end{align*}
and 
\begin{align*}
2(Q_1 w, Q_2 w)_{L^2(Q)} &= -\int_Q 4s (\na\al\!\cdot\!\na w_{{\rm im}}) (\pa_t w_{{\rm re}}) dxdt - \int_{Q} 2s(\De\al) w_{{\rm im}}(\pa_t w_{{\rm re}}) dxdt - \int_{Q} 4s(\na\al\!\cdot\!\na w_{{\rm im}})\De w_{{\rm im}} dxdt \\
&- \int_{Q} 2s (\De\al) w_{{\rm im}} \De w_{{\rm im}} dxdt - \int_{Q} 4s^3 |\na\al|^2(\na\al\!\cdot\!\na w_{{\rm im}})w_{{\rm im}} dxdt - \int_{Q} 2s^3 |\na\al|^2(\De\al) |w_{{\rm im}}|^2 dxdt, \\
&=: \sum_{k=1}^6 J_k,
\end{align*}
hence are left with the task of computing $I_k$ and $J_k$ for $j=1,\ldots,6$. We proceed by integration by parts and find through direct calculations that
\begin{align*}
I_1 &= \int_Q 4s (\na\al\!\cdot\!\na w_{{\rm re}}) (\pa_t w_{{\rm im}}) dxdt\\
&=-\int_Q 4s(\De\al) w_{{\rm re}}(\pa_t w_{{\rm im}}) dxdt - \int_Q 4s (\na\al\!\cdot\!\na (\pa_t w_{{\rm im}}))w_{{\rm re}} dxdt\\
&=-\int_{\Om} 4s(\na\al\!\cdot\!\na w_{{\rm im}})w_{{\rm re}} dx\Big|_{t=0}^{t=T} + \int_Q 4s (\na(\pa_t\al)\na w_{{\rm im}})w_{{\rm re}} dxdt + \int_Q 4s (\na\al\!\cdot\!\na w_{{\rm im}})(\pa_t w_{{\rm re}}) dxdt \\
&-\int_Q 4s(\De\al) w_{{\rm re}}(\pa_t w_{{\rm im}}) dxdt\\
&= -\int_{\Om} 2s(\na\al\!\cdot\!\na w_{{\rm im}})w_{{\rm re}} dx\Big|_{t=0}^{t=T} +\int_{\Om} 2s(\na\al\!\cdot\!\na w_{{\rm re}})w_{{\rm im}} dx\Big|_{t=0}^{t=T} + \int_{\Om} 2s(\De\al) w_{{\rm im}}w_{{\rm re}} dx\Big|_{t=0}^{t=T} \\
&+ \int_Q 4s (\na(\pa_t\al)\na w_{{\rm im}})w_{{\rm re}} dxdt + \int_Q 4s (\na\al\!\cdot\!\na w_{{\rm im}})(\pa_t w_{{\rm re}}) dxdt -\int_Q 4s(\De\al) w_{{\rm re}}(\pa_t w_{{\rm im}}) dxdt, \\
\\
I_2 &= \int_{Q} 2s(\De\al) w_{{\rm re}}(\pa_t w_{{\rm im}}) dxdt, \\
\\
I_3 &= - \int_{Q} 4s(\na\al\!\cdot\!\na w_{{\rm re}})\De w_{{\rm re}} dxdt \\
&=\int_\Sigma -4s \na\al\!\cdot\!\na w_{{\rm re}} \f{\pa w_{{\rm re}}}{\pa\nu} d\Sigma + \int_Q 4s\sum_{i,j=1}^n(\pa_i\pa_j\al)(\pa_i w_{{\rm re}})(\pa_j w_{{\rm re}}) dxdt + \int_Q 2s \na\al\!\cdot\!\na |\na w_{{\rm re}}|^2 dxdt \\
&=\int_\Sigma \left( -4s \na\al\!\cdot\!\na w_{{\rm re}} \f{\pa w_{{\rm re}}}{\pa\nu} + 2s\f{\pa\al}{\pa\nu}|\na w_{{\rm re}}|^2\right) d\Sigma + \int_Q 4s\sum_{i,j=1}^n(\pa_i\pa_j\al)(\pa_i w_{{\rm re}})(\pa_j w_{{\rm re}}) dxdt - \int_Q 2s\De\al |\na w_{{\rm re}}|^2 dxdt \\
&= -\int_\Sigma 2s \f{\pa\al}{\pa\nu}\left|\f{\pa w_{{\rm re}}}{\pa\nu}\right|^2 d\Sigma + \int_Q 4s\sum_{i,j=1}^n(\pa_i\pa_j\al)(\pa_i w_{{\rm re}})(\pa_j w_{{\rm re}}) dxdt - \int_Q 2s\De\al |\na w_{{\rm re}}|^2 dxdt, \\
\\
I_4 &= - \int_{Q} 2s (\De\al) w_{{\rm re}} \De w_{{\rm re}} dxdt\\
&= \int_Q 2s(\De\al) |\na w_{{\rm re}}|^2 dxdt + \int_Q 2s\na(\De\al)\na w_{{\rm re}} w_{{\rm re}} dxdt \\
&= \int_Q 2s(\De\al) |\na w_{{\rm re}}|^2 dxdt - \int_Q s(\De^2\al) |w_{{\rm re}}|^2 dxdt,\\
\\
I_5 &= - \int_{Q} 4s^3 |\na\al|^2(\na\al\!\cdot\!\na w_{{\rm re}})w_{{\rm re}} dxdt\\
&= \int_{Q} 2s^3 (\nabla \cdot (|\na\al|^2\na\al)) |w_{{\rm re}}|^2 dxdt, \\
\\
I_6 &= - \int_{Q} 2s^3 |\na\al|^2(\De\al) |w_{{\rm re}}|^2 dxdt, \\
\\
J_1 &= -\int_Q 4s (\na\al\!\cdot\!\na w_{{\rm im}}) (\pa_t w_{{\rm re}}) dxdt, \\
\\
J_2 &= - \int_{Q} 2s(\De\al) w_{{\rm im}}(\pa_t w_{{\rm re}}) dxdt \\
&= - \int_\Om 2s(\De\al) w_{{\rm im}} w_{{\rm re}} dx\Big|_{t=0}^{t=T} + \int_Q 2s(\De\al)(\pa_t w_{{\rm im}})w_{{\rm re}} dxdt + \int_Q 2s \pa_t (\De\al) w_{{\rm im}} w_{{\rm re}} dxdt, \\
\\
J_3 &= - \int_{Q} 4s(\na\al\!\cdot\!\na w_{{\rm im}})\De w_{{\rm im}} dxdt \\
&=\int_\Sigma -4s \na\al\!\cdot\!\na w_{{\rm im}} \f{\pa w_{{\rm im}}}{\pa\nu} d\Sigma + \int_Q 4s\sum_{i,j=1}^n(\pa_i\pa_j\al)(\pa_i w_{{\rm im}})(\pa_j w_{{\rm im}}) dxdt + \int_Q 2s \na\al\!\cdot\!\na |\na w_{{\rm im}}|^2 dxdt \\
&= -\int_\Sigma 2s \f{\pa\al}{\pa\nu}\left|\f{\pa w_{{\rm im}}}{\pa\nu}\right|^2 d\Sigma + \int_Q 4s\sum_{i,j=1}^n(\pa_i\pa_j\al)(\pa_i w_{{\rm im}})(\pa_j w_{{\rm im}}) dxdt - \int_Q 2s\De\al |\na w_{{\rm im}}|^2 dxdt, \\
\\
J_4 &= - \int_{Q} 2s (\De\al) w_{{\rm im}} \De w_{{\rm im}} dxdt \\
&= \int_Q 2s (\De\al) |\na w_{{\rm im}}|^2 dxdt + \int_Q 2s \na(\De\al) \cdot \na w_{{\rm im}} w_{{\rm im}} dxdt \\
&= \int_Q 2s(\De\al) |\na w_{{\rm im}}|^2 dxdt - \int_Q s(\De^2\al) |w_{{\rm im}}|^2 dxdt,\\
\\
J_5 &= - \int_{Q} 4s^3 |\na\al|^2(\na\al\!\cdot\!\na w_{{\rm im}})w_{{\rm im}} dxdt \\
&= \int_{Q} 2s^3 ( \nabla \cdot (|\na\al|^2\na\al) ) |w_{{\rm im}}|^2 dxdt, \\
\\
J_6 &= - \int_{Q} 2s^3 |\na\al|^2(\De\al) |w_{{\rm im}}|^2 dxdt.
\end{align*}
Therefore we have
$$
2(P_1 w, P_2 w)_{L^2(Q)} + 2(Q_1 w, Q_2 w)_{L^2(Q)}= \sum_{k=1}^{6} (I_k + J_k) =: {\mathrm{Main}}_1 + {\mathrm{Main}}_2 + {\mathrm{Lower}} + {\mathrm{Bndry}},
$$
with
\begin{align*}
{\mathrm{Main}}_1 &:= \int_Q 4s\sum_{i,j=1}^n (\pa_i\pa_j\al) \big[(\pa_i w_{{\rm re}})(\pa_j w_{{\rm re}}) + (\pa_i w_{{\rm im}})(\pa_j w_{{\rm im}})\big] dxdt, \\
{\mathrm{Main}}_2 &:= \int_Q 2s^3\big( \na|\na\al|^2 \na\al \big) |w|^2 dxdt\\
&= \int_Q 2s^3\sum_{i,j=1}^n (\pa_i\pa_j\al)(\pa_i\al)(\pa_j\al) |w|^2 dxdt,\\
{\mathrm{Lower}} &:= \int_Q 4s (\na(\pa_t\al)\na w_{{\rm im}})w_{{\rm re}} dxdt - \int_Q s(\De^2\al) |w|^2 dxdt + \int_Q 2s \pa_t (\De\al) w_{{\rm im}} w_{{\rm re}} dxdt,\\
{\mathrm{Bndry}} &:= -\int_\Sigma 2s \f{\pa\al}{\pa\nu}\left(\left|\f{\pa w_{{\rm re}}}{\pa\nu}\right|^2 + \left|\f{\pa w_{{\rm im}}}{\pa\nu}\right|^2 \right) d\Sigma + \int_{\Om} 2s\left[(\na\al\!\cdot\!\na w_{{\rm re}})w_{{\rm im}} - (\na\al\!\cdot\!\na w_{{\rm im}})w_{{\rm re}}\right] dx\Big|_{t=0}^{t=T}.
\end{align*}

Let us now introduce the weight functions
\begin{align}
\label{con:weight}
\al(x,t) = \f{e^{\la\be(x)} -e^{\la K}}{l^2(t)}\ \mbox{and}\ \vp(x,t) = \f{e^{\la\be(x)}}{l^2(t)},
\end{align}
where $\be \in C^4(\ov{\Om})$ is nonnegative and has no critical point, i.e.
\begin{align}
\label{con:beta1}
\be(x)\ge 0,\ \quad \abs{\na\be(x)}\ge c_0>0,\ \forall x\in\Om,
\end{align}
where $K:=2\sup_{x\in\Om}\be(x)$ and $l\in C^1[0,T]$ is nonnegative, attains its maximum at the origin and vanishes at $T$, i.e.,
\begin{align}
\label{con:l}
l(T) = 0, \qquad l(0)> l(t)\ge 0, \quad \forall t\in (0,T].
\end{align}
We assume in addition that $\beta$ is pseudo-convex condition with respect to the Laplace operator, in the sense that there exist two constants $\la_1 \in (0,+\infty)$ and $\ep \in (0,+\infty)$ such that 
\begin{align}
\label{con:beta2}
\la \abs{\na \be\cdot \xi}^2 + D^2\be(\xi,\xi) \ge \ep \abs{\xi}^2,\ \xi\in \Bbb{R}^n,\ \la \in (\la_1,+\infty).
\end{align}
Next we define the observation zone where the Neumann data used by the analysis of the inverse problems examined in this text, are measured, as the sub-boundary 
$$
\Ga_0 := \{x\in\pa\Om:  \na \be(x) \cdot \nu(x) \ge 0\}.
$$
\begin{rem}
At this point it is worth mentioning that there exist actual functions $\be$ and $l$ fulfilling the conditions \eqref{con:beta1}, \eqref{con:l} and \eqref{con:beta2}. As a matter of fact, for any fixed $x_0\notin \ov{\Om}$, we may choose
$\be(x) := \abs{x-x_0}^2$ for all $x\in \ov{\Om}$ and $l(t) := (T+t)(T-t)$ for all $t \in [0,T]$.
In this case, the observation zone $\Ga_0$ coincides with the $x_0$-shadowed face of the boundary $\pa \Om$, i.e.
$\Ga_0 = \{x\in\pa\Om:\ (x-x_0)\cdot\nu(x)\ge 0\}$.
\end{rem}
From the very definition of $\al$, we see that $\lim_{t \to T} (\vp w)(\cdot,t)=0$, and for all  $i,j=1,\ldots,n$, that
\begin{align*}
&\na\al =\na\vp=\la\vp\na\be, \qquad \pa_i\al = \pa_i \vp = \la\vp(\pa_i \be), \qquad \pa_i\pa_j \al = \pa_i\pa_j\vp =\la^2\vp(\pa_i\be)(\pa_j\be) +\la\vp(\pa_i\pa_j \be)\\
&|\pa_t \al| =\left|-\f{2l^\prime(e^{\la\be}-e^{\la K})}{l^3}\right|\le C_{\la}\vp^{\f{3}{2}}, \qquad |\na(\pa_t\al)| = |\la(\pa_t\vp)\na \be| \le C_{\la}\vp^{\f{3}{2}},\qquad |\pa_t (\De\al)|\le C_{\la}\vp^{\f{3}{2}}.
\end{align*}
Here and henceforth, $C$ (resp., $C_\la$) denotes a generic constant that depends only on $\ep$, $c_0$ and $l(0)$ (resp., $\ep$, $c_0$, $\|\be\|_{L^\infty(\Om)}$, $l$ and $\la$). In any case, $C$ and $C_\la$ are independent of $s$. Therefore we have
\begin{align*}
{\mathrm{Main}}_1 &= \int_Q 4s\sum_{i,j=1}^n (\pa_i\pa_j\al) \big[(\pa_i w_{{\rm re}})(\pa_j w_{{\rm re}}) + (\pa_i w_{{\rm im}})(\pa_j w_{{\rm im}})\big] dxdt\\
&= \int_Q 4s\la\vp\left[\la|\na\be\!\cdot\!\na w_{{\rm re}}|^2 + \la|\na\be\!\cdot\!\na w_{{\rm im}}|^2 + D^2\be (\na w_{{\rm re}},\na w_{{\rm re}}) + D^2\be (\na w_{{\rm im}},\na w_{{\rm im}})\right] dxdt \\
&\ge 4\ep\int_Q s\la\vp |\na w|^2 dxdt,\\
{\mathrm{Main}}_2 &= \int_Q 2s^3\sum_{i,j=1}^n (\pa_i\pa_j\al)(\pa_i\al)(\pa_j\al) |w|^2 dxdt\\
&=\int_Q 2s^3\la^3\vp^3\left[\la |\na\be|^4 + D^2\be(\na\be,\na\be) \right] |w|^2 dxdt\\
&\ge 2\ep \int_Q 2s^3\la^3\vp^3 |\na\be|^2|w|^2 dxdt\\
&\ge 2\ep c_0^2 \int_Q 2s^3\la^3\vp^3 |w|^2 dxdt,\\
|{\mathrm{Lower}}| &= \left|\int_Q 4s (\na(\pa_t\al)\na w_{{\rm im}})w_{{\rm re}} dxdt - \int_Q s(\De^2\al) |w|^2 dxdt + \int_Q 2s \pa_t (\De\al) w_{{\rm im}} w_{{\rm re}} dxdt\right|\\
&\le C_{\la}s^{-\f{1}{2}}\int_Q (s\vp |\na w|^2 + s^3\vp^3 |w|^2) dxdt,\\
{\mathrm{Bndry}} &= -\int_\Sigma 2s \f{\pa\al}{\pa\nu}\left(\left|\f{\pa w_{{\rm re}}}{\pa\nu}\right|^2 + \left|\f{\pa w_{{\rm im}}}{\pa\nu}\right|^2 \right) d\Sigma + \int_{\Om} 2s\left[(\na\al\!\cdot\!\na w_{{\rm re}})w_{{\rm im}} - (\na\al\!\cdot\!\na w_{{\rm im}})w_{{\rm re}}\right] dx\Big|_{t=0}^{t=T}\\
&= -\int_\Sigma 2s\la\vp \f{\pa\be}{\pa\nu}\left|\f{\pa w}{\pa\nu}\right|^2 d\Sigma - \int_{\Om} 2s\la\vp\left[(\na\be\!\cdot\!\na u_{{\rm re}})u_{{\rm im}} - (\na\be\!\cdot\!\na u_{{\rm im}})u_{{\rm re}}\right]e^{2s\al} dx\Big|_{t=0}\\
&\ge -\int_{\Sigma_0} 2s\la\vp \f{\pa\be}{\pa\nu}\left|\f{\pa w}{\pa\nu}\right|^2 d\Sigma - \int_{\Om} 2s\la\vp\left[(\na\be\!\cdot\!\na u_{{\rm re}})u_{{\rm im}} - (\na\be\!\cdot\!\na u_{{\rm im}})u_{{\rm re}}\right]e^{2s\al} dx\Big|_{t=0},
\end{align*}
where $\Sigma_0 := (0,T) \times \Gamma_0$. This and \eqref{eq1} imply
\begin{align*}
&\|R_1 w\|_{L^2(Q)}^2 + \|R_2 w\|_{L^2(Q)}^2 + \|s^{\f{1}{2}}\vp^{\f{1}{2}}\na w\|_{L^2(Q)}^2 + \|s^{\f{3}{2}}\vp^{\f{3}{2}} w\|_{L^2(Q)}^2 \le C\|R_3 w\|_{L^2(Q)}^2 + C_{\la}s\left\|\vp^{\f{1}{2}}\f{\pa w}{\pa\nu}\right\|_{L^2(\Sigma_0)}^2\\
&+ C_{\la}s^{-\f{1}{2}}\left(\|s^{\f{1}{2}}\vp^{\f{1}{2}}\na w\|_{L^2(Q)}^2 + \|s^{\f{3}{2}}\vp^{\f{3}{2}} w\|_{L^2(Q)}^2\right)  + C\left|\int_{\Om} 2s\la\vp\left[(\na\be\!\cdot\!\na u_{{\rm re}})u_{{\rm im}} - (\na\be\!\cdot\!\na u_{{\rm im}})u_{{\rm re}}\right]e^{2s\al} dx\Big|_{t=0}\right|
\end{align*}
for all $\la\ge \la_2 :=\max\{1,\la_1\}$ and all $s \in (0,+\infty)$. Further, bearing in mind that
$R_3 w = R w - is(\pa_t \al)w$ and $R w = e^{s\al} Lu$, we infer from the above inequality that
\begin{align*}
&\|R_1 w\|_{L^2(Q)}^2 + \|R_2 w\|_{L^2(Q)}^2 + \|s^{\f{1}{2}}\vp^{\f{1}{2}}\na w\|_{L^2(Q)}^2 + \|s^{\f{3}{2}}\vp^{\f{3}{2}} w\|_{L^2(Q)}^2 \le C\|R w\|_{L^2(Q)}^2 + C_{\la}s\left\|\vp^{\f{1}{2}}\f{\pa w}{\pa\nu}\right\|_{L^2(\Sigma_0)}^2\\
&+ C_{\la}s^{-\f{1}{2}}\left(\|s^{\f{1}{2}}\vp^{\f{1}{2}}\na w\|_{L^2(Q)}^2 + \|s^{\f{3}{2}}\vp^{\f{3}{2}} w\|_{L^2(Q)}^2\right)  + C\int_{\Om} 2s\la\vp\left|(\na\be\!\cdot\!\na u_{{\rm re}})u_{{\rm im}} - (\na\be\!\cdot\!\na u_{{\rm im}})u_{{\rm re}}\right|e^{2s\al} dx\Big|_{t=0}.
\end{align*}
Thus, going back to $u=e^{-s \al } w$ and taking $\la\ge \la_3:=\max\{\la_2, 2(\ln{2})K^{-1}\}$ and $s\ge s_1(\la):= 4C_{\la}^2>0$ in such a way that the low order term $C_{\la}s^{-\f{1}{2}}\left(\|s^{\f{1}{2}}\vp^{\f{1}{2}}\na w\|_{L^2(Q)}^2 + \|s^{\f{3}{2}}\vp^{\f{3}{2}} w\|_{L^2(Q)}^2\right)$ in the right-hand side of the above inequality is absorbed by $\|s^{\f{1}{2}}\vp^{\f{1}{2}}\na w\|_{L^2(Q)}^2 + \|s^{\f{3}{2}}\vp^{\f{3}{2}} w\|_{L^2(Q)}^2$ in the left-hand side, we get for all $s \in [s_1,+\infty)$ that
\begin{align}
&\|R_1 (ue^{s\al})\|_{L^2(Q)}^2 + \|R_2 (ue^{s\al})\|_{L^2(Q)}^2 + \|s^{\f{1}{2}}\vp^{\f{1}{2}}e^{s\al}\na u\|_{L^2(Q)}^2 + \|s^{\f{3}{2}}\vp^{\f{3}{2}}e^{s\al} u\|_{L^2(Q)}^2 \le C\|e^{s\al}Lu\|_{L^2(Q)}^2 \nonumber\\
&\hspace{2.5cm} + C_{\la}s\left\|\vp^{\f{1}{2}}e^{s\al}\f{\pa u}{\pa\nu}\right\|_{L^2(\Sigma_0)}^2 + C_{\la}\int_{\Om} 2s\vp\left|(\na\be\!\cdot\!\na u_{{\rm re}})u_{{\rm im}} - (\na\be\!\cdot\!\na u_{{\rm im}})u_{{\rm re}}\right|e^{2s\al} dx\Big|_{t=0}.
\label{c10}
\end{align}
Further, for all $\lambda \in [\lambda_3,+\infty)$ and all $(x,t) \in Q$, we notice that $s \vp(x,t)e^{2s\al(x,t)} \le s l^{-2}(t) e^{\la K} e^{-2s l^{-2}(t)e^{\la K}}$, in such a way that
\begin{equation}
\label{c11}
s\left\|\vp^{\f{1}{2}}e^{s\al}\f{\pa u}{\pa\nu}\right\|_{L^2(\Sigma_0)}^2 \le C_{\la}\left\|\f{\pa u}{\pa\nu}\right\|_{L^2(\Sigma_0)}^2.
\end{equation}
Thus, taking into account that 
$\vp(x,0)= e^{\la\be(x)}l^{-2}(0)\le e^{\f{\la K}{2}}l^{-2}(0)$ and $\vp(x,t)=e^{\la\be(x)}l^{-2}(t)\ge l^{-2}(0)$ for all $x \in \Om$ and all $t \in [0,T]$, we infer from \eqref{c10}-\eqref{c11} that
\begin{align*}
&\|R_1 (ue^{s\al})\|_{L^2(Q)}^2 + \|R_2 (ue^{s\al})\|_{L^2(Q)}^2 + s\|e^{s\al}\na u\|_{L^2(Q)}^2 + s^3\|e^{s\al} u\|_{L^2(Q)}^2\\
&\le C\|e^{s\al}Lu\|_{L^2(Q)}^2 + C_{\la}\left\|\f{\pa u}{\pa\nu}\right\|_{L^2(\Sigma_0)}^2 + C_{\la}\int_{\Om} 2s\left|(\na\be\!\cdot\!\na u_{{\rm re}})u_{{\rm im}} - (\na\be\!\cdot\!\na u_{{\rm im}})u_{{\rm re}}\right|e^{2s\al} dx\Big|_{t=0}\\
&\le C\|e^{s\al}Lu\|_{L^2(Q)}^2 + C_{\la}\left\|\f{\pa u}{\pa\nu}\right\|_{L^2(\Sigma_0)}^2 + C_{\la}s Z(u_0),\ s \in [s_1,+\infty),
\end{align*}
with
\begin{equation}
\label{defZ}
Z(u_0) := \int_{\Om} e^{2s\al(x,0)} \abs{\left( \ov{u_0} \na\be\!\cdot\!\na u_0 - u_0 \na\be\!\cdot\!\na \ov{u_0} \right)(x)}  dx.
\end{equation}

Thus we have proved the following:

\begin{thm}
\label{thm:CE}
Let $\al$ be defined by \eqref{con:weight} where $\be\in C^4(\ov{\Om})$ and $l\in C^1[0,T]$ fulfill the conditions \eqref{con:beta1}-\eqref{con:beta2} for some fixed $\la \in (0,+\infty)$. Then there exist two positive constants $s_0$ and $C$, both of them depending only on $\ep$, $c_0$, $\la$, $l(0)$, $\|\be\|_{L^\infty(\Om)}$ and $\|l^\prime\|_{L^\infty(\Om)}$,
such that the estimate
\begin{align}
\nonumber
&\|R_1 (ue^{s\al})\|_{L^2(Q)}^2 + \|R_2 (ue^{s\al})\|_{L^2(Q)}^2 + s\|e^{s\al}\na u\|_{L^2(Q)}^2 + s^3\|e^{s\al} u\|_{L^2(Q)}^2\\
\label{CE}
&\hspace{2cm}\le C\|e^{s\al}Lu\|_{L^2(Q)}^2 + C\left\|\f{\pa u}{\pa\nu}\right\|_{L^2(\Sigma_0)}^2 + Cs\int_{\Om} i\left[(\na\be\!\cdot\!\na u)\ov{u} - (\na\be\!\cdot\!\na \ov{u})u \right]e^{2s\al} dx\Big|_{t=0}
\end{align}
holds for all $s \in [s_0,+\infty)$ and any function $u\in L^2(0,T; H_0^1(\Om))$ satisfying $Lu\in L^2(\Om\times (0,T))$ and $\pa_\nu u\in L^2(0,T; L^2(\pa\Om))$.  Here the operators $R_1$, $R_2$ and $L$ are defined in \eqref{equ-gov}-\eqref{def:R}.  
\end{thm}

\section{Proof of Theorems \ref{thm:stb}, \ref{thm:stb2} and \ref{thm:stb3}}
\label{sec-pr}

\subsection{Preliminary estimate}
Let us recall that $u_j$, $j=1,2$, is the $H^2(0,T;L^2(\Omega)) \cap H^1(0,T;H^2(\Omega))$-solution to \eqref{sys:ori} where $(A,\rho)$ is replaced by $(A_j,\rho_j)$. Thus, taking the difference of the two systems, we get that $u:=u_1-u_2$ solves
\begin{equation*}
\left\{
\begin{aligned}
& (-i\pa_t - \De_{A_1} + \rho_1) u = 2i A\!\cdot\!\na u_2 -\left(\rho + (A_1+A_2)\!\cdot\! A -i \na\!\cdot\! A \right) u_2 && \text{in}\ Q \\
& u= 0 \qquad && \text{on}\ \Sigma \\
& u(\cdot,0) = 0 \qquad && \text{on}\ \Om,
\end{aligned}
\right.
\end{equation*}
with $\rho:=\rho_1-\rho_2$ and $A:=A_1-A_2$. Further, differentiating the above system w.r.t. the time variable $t$, yields
\begin{equation}
\label{sys:maglinv}
\left\{
\begin{aligned}
& (-i\pa_t  - \De_{A_1}  + \rho_1 ) v = 2i A\!\cdot\!\na \pa_t u_2 -\left(\rho + (A_1+A_2)\!\cdot\! A - i \na\!\cdot\! A \right) \pa_t u_2 \qquad && \text{in}\ Q \\
& v = 0 \qquad && \text{on}\ \Sigma \\
& v(\cdot,0) = -2A \!\cdot\!\na u_0 - i\left(\rho + (A_1+A_2)\!\cdot\! A - i \na\!\cdot\! A \right) u_0 \qquad && \text{in}\ \Om
\end{aligned}
\right.
\end{equation}
with $v:=\pa_t u$. All the above computations make sense as we have $u \in H^2(0,T; L^2(\Om))\cap H^1(0,T; H^2(\Om) \cap H_0^1(\Om))$, and hence
\begin{equation}
\label{con:reg2}
v\in H^1(0,T; L^2(\Om))\cap L^2(0,T; H^2(\Om) \cap H_0^1(\Om)).
\end{equation}
Moreover, we have the following estimate:
\begin{equation}
\label{eq:lem}
\norm{e^{s \al(\cdot,0)} v (\cdot,0)}_{L^2(\Om)}^2 \le s^{-\f{3}{2}} \left(\norm{R_1 e^{s\al} v}_{L^2(Q)}^2 + s^3 \norm{e^{s \al} v}_{L^2(Q)}^2 \right),
\ s \in (0,+\infty),
\end{equation}
where $R_1$ is defined by \eqref{def:R}.
This can be seen through direct calculations. Indeed, in light of \eqref{con:weight}-\eqref{con:l} we see that $\lim_{t \stackrel{<}{\to}T} \al(x,t)=-\infty$ for all $x \in \Om$, whence $\lim_{t \stackrel{<}{\to}T} e^{s \alpha(\cdot,t)} v(\cdot,t)=0$ in $L^2(\Om)$. As a consequence we have
$$
\norm{e^{s \alpha(\cdot,0)} v(\cdot,0)}_{L^2(\Om)}^2 = - \int_0^T \f{d}{dt} \norm{e^{s \al(\cdot,t)} v(\cdot,t)}_{L^2(\Om)}^2 = 
-\int_Q \left(  w \ov{\pa_t w} + (\pa_t w) \ov{w} \right) dx\ dt,
$$
with $w:=\e^{s \alpha} v$. Further, as $\pa_t w = i \left( R_1w  + \De w + s^2 \abs{\na\al}^2 w \right)$ from the very definition of $R_1$, it holds true that\begin{align*}
\norm{e^{s \alpha(\cdot,0)} v(\cdot,0)}_{L^2(\Om)}^2 &=  i \int_Q \left( w( \ov{R_1 w} + \De \ov{w} + s^2 \abs{\na \al}^2 \ov{w}) - (R_1 w + \De w + s^2 \abs{\na\al}^2 w) \ov{w} \right) dx\ dt \\
& = i \int_Q \left( w \ov{R_1w} - (R_1w)\ov{w} + w \De \ov{w} - \ov{w} \De w \right) dx\ dt
\end{align*}
Finally, bearing in mind that $\int_Q \left(  w \De \ov{w} - \ov{w} \De w \right) dx\ dt =0$, we end up getting
$$
\norm{e^{s \alpha(\cdot,0)} v(\cdot,0)}_{L^2(\Om)}^2 =
2 \mathrm{Re} \langle s^{\f{3}{4}}iw, s^{-\f{3}{4}}R_1w \rangle_{L^2(Q)} \le s^{-\f{3}{2}} \norm{R_1 w}_{L^2(Q)}^2 + s^{\f{3}{2}} \norm{w}_{L^2(Q)}^2,
$$
with the help of the Cauchy-Schwarz inequality. This leads to \eqref{eq:lem}.

\subsection{Proof of Theorem \ref{thm:stb}}
Let us rewrite \eqref{sys:maglinv} in the context of Theorem \ref{thm:stb}, where $A_1=A_2=0$. We get
\begin{equation}
\label{sys:linv}
\left\{
\begin{aligned}
& (-i\pa_t  - \De  + \rho_1) v = -\rho \pa_t u_2 \qquad && \text{in}\ Q \\
& v = 0 \qquad && \text{on}\ \Sigma \\
& v(\cdot,0) = -i \rho u_0 \qquad && \text{in}\ \Om.
\end{aligned}
\right.
\end{equation}
With reference to \eqref{equ-gov}, the first line of \eqref{sys:linv} reads $Lv= -\rho \pa_t u_2 - \rho_1 v$. Thus we have $Lv\in L^2(0,T;L^2(\Om))$, with $v\in L^2(0,T;H_0^1(\Om))$ and $\pa_\nu v\in L^2(0,T; L^2(\pa\Om))$, according to \eqref{con:reg2}. Therefore, Theorem \ref{thm:CE} yields
\begin{align*}
&\norm{R_1 (e^{s\al}v)}_{L^2(Q)}^2 + \norm{R_2 (e^{s\al}v)}_{L^2(Q)}^2 + s \norm{e^{s\al}\na v}_{L^2(Q)}^2 + s^3 \norm{e^{s\al} v}_{L^2(Q)}^2 \\ 
&\le C \left( \norm{\pa_\nu v}_{L^2(\Sigma_0)}^2 + \norm{e^{s\al}(\rho \pa_t u_2 + \rho_1 v)}_{L^2(Q)}^2 +  s Z(-i \rho u_0)   \right) \\
& \le C \left(\norm{\pa_\nu v}_{L^2(\Sigma_0)}^2 + \norm{e^{s\al}\rho}_{L^2(Q)}^2 + \norm{e^{s\al}v}_{L^2(Q)}^2 + s Z(-i \rho u_0) \right),\ s \in (s_0,+\infty),
\end{align*}
where $Z$ is defined in \eqref{defZ}. In the last line we used the energy estimate \eqref{c-es} where $u_2$ is substituted for $u$. 

Further, upon taking $s_1 \in [\max(s_0,1),+\infty)$ so large that $s_1^3 \geq 2C$ (in such a way that $C \norm{e^{s_1\al}v}_{L^2(Q)}^2$ is absorbed by $s_1^3 \norm{e^{s_1\al} v}_{L^2(Q)}^2$), we get that
\begin{align*}
& \norm{R_1 (e^{s\al}v)}_{L^2(Q)}^2 + \norm{R_2 (e^{s\al}v)}_{L^2(Q)}^2 + s\norm{e^{s\al}\na v}_{L^2(Q)}^2 + s^3\norm{e^{s\al} v}_{L^2(Q)}^2\\
& \le C \left( \norm{\pa_\nu v}_{L^2(\Sigma_0)}^2 + \norm{e^{s\al}\rho}_{L^2(Q)}^2 + s Z(-i \rho u_0) \right),\ s \in (s_1,+\infty).
\end{align*}
From this and \eqref{eq:lem}, it then follows that
\begin{align*}
\norm{e^{s\al(\cdot,0)}v(\cdot,0)}_{L^2(\Om)}^2 \le C\left(\norm{\pa_\nu v}_{L^2(\Sigma_0)}^2 + s^{-\f{3}{2}} \norm{e^{s\al}\rho}_{L^2(Q)}^2 + s^{-\f{1}{2}} Z(-i \rho u_0) \right),\ s \in (s_1,+\infty).
\end{align*}
Moreover, we have
$\norm{e^{s\al(\cdot,0)}v(\cdot,0)}_{L^2(\Om)}^2 =  \norm{e^{s\al(\cdot,0)} \rho u_0}_{L^2(\Om)}^2 \ge r_0^2 \norm{e^{s\al(\cdot,0)} \rho}_{L^2(\Om)}^2$ for all $s \in \R$, from \eqref{ass:u_0} and the third line of \eqref{sys:linv}, hence
\begin{align}
\label{eq:inv1}
\norm{e^{s\al(\cdot,0)}\rho}_{L^2(\Om)}^2 \le C\left(\norm{\pa_\nu v}_{L^2(\Sigma_0)}^2 +s^{-\f{3}{2}} \norm{e^{s\al}\rho}_{L^2(Q)}^2 + s^{-\f{1}{2}} Z(-i \rho u_0) \right),\ s \in (s_1,+\infty).
\end{align}
We are left with the task of estimating the two last terms appearing in the right-hand side of \eqref{eq:inv1}. For the first one, we take advantage of the fact, arising from \eqref{con:weight} and \eqref{con:l}, that 
\begin{equation}
\label{c-weight}
\alpha(x,t) \le \alpha(x,0),\ (x,t) \in Q.
\end{equation}
This yields
\begin{equation}
\label{eq:inv2}
\norm{e^{s\al}\rho}_{L^2(Q)}^2 \le \norm{e^{s\al(\cdot,0)}\rho}_{L^2(Q)}^2 = T\norm{e^{s\al(\cdot,0)}\rho}_{L^2(\Om)}^2,\ s \in [0,+\infty).
\end{equation}
For the second term, we infer from \eqref{defZ} and the last line of \eqref{sys:linv} that
\begin{align*}
Z(-i \rho u_0)=
& \int_{\Om} e^{2s\al(x,0)} \left( \abs{u_0}^2 (\ov{\rho}\na\be\!\cdot\!\na\rho - \rho\na\be\!\cdot\!\na\ov{\rho}) + \abs{\rho}^2 (\ov{u_0}\na\be\!\cdot\!\na u_0 - u_0\na\be\!\cdot\!\na\ov{u_0}) \right)(x)  dx\\
\hspace{0cm}=& \int_{\Om} e^{2s\al(x,0)}\abs{\na\be\cdot\left( \abs{u_0}^2(\ov{\rho}\na\rho - \rho\na\ov{\rho}) + \abs{\rho}^2(\ov{u_0}\na u_0 - u_0\na\ov{u_0}) \right)(x)}  dx.
\end{align*}
Thus, with reference to \eqref{c-el}, entailing 
$\abs{\na\be\cdot(\ov{\rho}\na\rho - \rho\na\ov{\rho})} \le C \abs{\rho}^2$ a.e. in $\Om$, we get
\begin{equation*}
Z(-i \rho u_0) \le C \norm{e^{s\al(\cdot,0)}\rho}_{L^2(\Om)}^2,\ s \in \R.
\end{equation*}
Now, putting this together with  \eqref{eq:inv1}-\eqref{eq:inv2}, we obtain
$$
\norm{e^{s\al(\cdot,0)}\rho}_{L^2(\Om)}^2
\le C \left( \norm{\pa_\nu v}_{L^2(\Sigma_0)}^2 +  (s^{-\f{1}{2}} + s^{-\f{3}{2}})\norm{e^{s\al(\cdot,0)}\rho}_{L^2(\Om)}^2 \right),\ s \in (s_1,+\infty).
$$ 
Thus, taking $s_2 \in (s_1,+\infty)$ so large that $s_2^{-\f{1}{2}} + s_2^{-\f{3}{2}}$ is not greater than, say, $\f{1}{2C}$, the above estimate immediately yields
\begin{equation}
\label{eq:inv3}
\norm{e^{s\al(\cdot,0)}\rho}_{L^2(\Om)}^2 \le C \norm{\pa_\nu v}_{L^2(\Sigma_0)}^2,\ s \in (s_2,+\infty).
\end{equation}
Finally, the desired result follows readily from this and the following estimate, arising from \eqref{con:weight}--\eqref{con:l}:
\begin{equation}
\label{e:0}
e^{s\al(x,0)}=e^{sl^{-2}(0)(e^{\la\be(x)}-e^{\la K})}\ge e^{sl^{-2}(0)(1-e^{\la K})} \in (0,+\infty),\ x \in \Om,\ s \in [0,+\infty). 
\end{equation}

\subsection{Proof of Theorem \ref{thm:stb2}}
\label{sec-stb2}
With reference to \eqref{def-deltaA}, \eqref{equ-gov} and the first line of \eqref{sys:maglinv}, we have
$$
Lv= 2i A\!\cdot\!\na \pa_t u_2 -\left(\rho + (A_1+A_2)\!\cdot\! A -i \na\!\cdot\! A \right) \pa_t u_2 + 2i A_1\!\cdot\!\na v - (\rho_1 + \abs{A_1}^2 - i \na\!\cdot\! A_1) v.
$$
Thus $v\in L^2(0,T;H_0^1(\Om))$, $L v\in L^2(0,T;L^2(\Om))$ and $\pa_\nu v\in L^2(0,T; L^2(\pa\Om))$, by\eqref{con:reg2}, and we get
$$
\norm{e^{s \al} L v}_{L^2(Q)}^2 \le C \left( \norm{e^{s\al}\rho}_{L^2(Q)}^2 + 
\norm{e^{s \al} A}_{L^2(Q)}^2 + \norm{e^{s \al} \na \cdot A}_{L^2(Q)}^2 + \norm{e^{s \al} v}_{L^2(Q)}^2 + \norm{e^{s \al} \na v}_{L^2(Q)}^2 \right),\ s \in \R,
$$
upon substituting $u_2$ for $u$ in \eqref{c-es}. Therefore, in light of the Carleman estimate of Theorem \ref{thm:CE},
we obtain for all $s \in (s_0,+\infty)$:
\begin{align*}
&  \norm{R_1 (e^{s\al} v)}_{L^2(Q)}^2 +  \norm{R_2 (e^{s\al} v)}_{L^2(Q)}^2 + s \norm{e^{s\al} \na v}_{L^2(Q)}^2 + s^3 \norm{e^{s\al} v}_{L^2(Q)}^2\\ 
& \le C \left( \norm{\pa_\nu v}_{L^2(\Sigma_0)}^2 + \norm{e^{s\al}\rho}_{L^2(Q)}^2 + 
\norm{e^{s \al} A}_{L^2(Q)}^2 + \norm{e^{s \al} \na \cdot A}_{L^2(Q)}^2 + \norm{e^{s \al} v}_{L^2(Q)}^2 + \norm{e^{s \al} \na v}_{L^2(Q)}^2 + s Z(v(\cdot,0)) \right).
\end{align*}
Here $Z(v(\cdot,0))$ is as in \eqref{defZ} with
$v(\cdot,0)=-2A \!\cdot\!\na u_0 - i\left(\rho + (A_1+A_2)\!\cdot\! A - i \na\!\cdot\! A \right) u_0$, from the third line of \eqref{sys:maglinv}. Next, taking $s_1 \in [s_0,+\infty)$ so large that $\min(s_1,s_1^3) \ge 2C$ then yields
\begin{align*}
& \norm{R_1 (e^{s\al} v)}_{L^2(Q)}^2 +  \norm{R_2 (e^{s\al} v)}_{L^2(Q)}^2 + s \norm{e^{s\al} \na v}_{L^2(Q)}^2 + s^3 \norm{e^{s\al} v}_{L^2(Q)}^2\\
&  \le C \left( \norm{\pa_\nu v}_{L^2(\Sigma_0)}^2 +  \norm{e^{s\al}\rho}_{L^2(Q)}^2 + \norm{e^{s\al}A}_{L^2(Q)}^2 + \norm{e^{s\al} \na\!\cdot\! A}_{L^2(Q)}^2 + s Z(v(\cdot,0)) \right),\ s \in (s_1,+\infty).
\end{align*}
This and \eqref{eq:lem} imply for all $s \in (s_1,+\infty)$,
\begin{equation}
\label{eq:inv21}
\norm{e^{s\al(\cdot,0)} v(\cdot,0)}_{L^2(\Om)}^2  \le Cs^{-\f{3}{2}} \left( \norm{\pa_\nu v}_{L^2(\Sigma_0)}^2 + \norm{e^{s\al}\rho}_{L^2(Q)}^2 + \norm{e^{s\al}A}_{L^2(Q)}^2 + \norm{e^{s\al} \na\!\cdot\! A }_{L^2(Q)}^2 + s Z(v(\cdot,0)) \right), 
\end{equation}
with
\begin{eqnarray}
& & \norm{e^{s\al}\rho}_{L^2(Q)}^2 + \norm{e^{s\al}A}_{L^2(Q)}^2 + \norm{e^{s\al} \na\!\cdot\! A}_{L^2(Q)}^2 \nonumber \\
&  \le & T \left( \norm{e^{s\al(\cdot,0)}\rho}_{L^2(\Om)}^2 + \norm{e^{s\al(\cdot,0)}A}_{L^2(\Om)}^2 + \norm{e^{s\al(\cdot,0)} \na\!\cdot\! A}_{L^2(\Om)}^2 \right),\ s \in (0,+\infty), \label{e:1}
\end{eqnarray}
in virtue of \eqref{c-weight}. Similarly, we have
\begin{equation} 
\label{eq:Zes}
Z(v(\cdot,0)) \le C\left(\norm{e^{s\al(\cdot,0)}\rho}_{L^2(\Om)}^2 + \norm{e^{s\al(\cdot,0)}A}_{L^2(\Om)}^2 + \norm{e^{s\al(\cdot,0)} \na\!\cdot\! A}_{L^2(\Om)}^2 \right),\ s \in \R.
\end{equation}
This can be seen from the third line of \eqref{sys:maglinv}, entailing
\begin{equation}
\label{e:2}
\norm{e^{s \alpha(\cdot,0)} v(\cdot,0)}_{L^2(\Om)} \le C \left( \norm{e^{s \alpha(\cdot,0)} \rho}_{L^2(\Omega)} + \norm{e^{s \alpha(\cdot,0)} A}_{L^2(\Om)} + \norm{e^{s \alpha(\cdot,0)} \nabla \cdot A}_{L^2(\Om)} \right),\ s \in \R,
\end{equation}
and from the identity $\nabla v(\cdot,0) = -(w_0 + z_0)$ with $w_0:= 2 \mathbb{D}^2_{u_0} A + (\nabla \cdot A) \nabla u_0 + i \left( \rho+(A_1+A_2) \cdot A \right) \nabla u_0 + i \mathbb{J}_{A_1+A_2} A u_0 $ and $z_0:=2 \mathbb{J}_{A} \nabla u_0 + u_0 \nabla (\nabla \cdot A) + i \left( \nabla \rho+ \mathbb{J}_{A} (A_1+A_2) \right) u_0$, where we have set $\mathbb{D}^2_{u_0}:=\left( \partial_{i j}^2 u_0 \right)_{1 \le i, j \le n}$ and $\mathbb{J}_{Y}:=\left( \partial_i y_j \right)_{1 \le i, j \le n}$ for all $Y=(y_j)_{1 \le j \le n} \in H^1(\Omega,\R^n)$.
As a matter of fact we have 
\begin{equation}
\label{e:2b}
\norm{e^{s \alpha(\cdot,0)} w_0}_{L^2(\Om)} \le C \left( \norm{e^{s \alpha(\cdot,0)} \rho}_{L^2(\Omega)} + \norm{e^{s \alpha(\cdot,0)} A}_{L^2(\Om)} + \norm{e^{s \alpha(\cdot,0)} \nabla \cdot A}_{L^2(\Om)} \right),\ s \in \R,
\end{equation}
and we infer from assumptions \eqref{h-em-a},  \eqref{h-em-b} and \eqref{h-em-c} that
\begin{align*}
\norm{e^{s \alpha(\cdot,0)} z_0}_{L^2(\Om)} & \le C \left( \norm{e^{s \alpha(\cdot,0)} \nabla \rho}_{L^2(\Omega)} + \norm{e^{s \alpha(\cdot,0)} \mathbb{J}_{A}}_{L^2(\Om)} + \norm{e^{s \alpha(\cdot,0)} \nabla (\nabla \cdot A)}_{L^2(\Om)} \right) \\
& \le C \left( \norm{e^{s \alpha(\cdot,0)} \rho}_{L^2(\Omega)} + \norm{e^{s \alpha(\cdot,0)} A}_{L^2(\Om)} + \norm{e^{s \alpha(\cdot,0)} \nabla \cdot A}_{L^2(\Om)} \right),\ s \in \R,
\end{align*}
so we get 
$$
\norm{e^{s \alpha(\cdot,0)} \nabla v(\cdot,0)}_{L^2(\Om)} \le C \left( \norm{e^{s \alpha(\cdot,0)} \rho}_{L^2(\Omega)} + \norm{e^{s \alpha(\cdot,0)} A}_{L^2(\Om)} + \norm{e^{s \alpha(\cdot,0)} \nabla \cdot A}_{L^2(\Om)} \right),\ s \in \R.
$$
This and \eqref{e:2} entail \eqref{eq:Zes} through the Cauchy-Schwarz inequality. 
Next, putting \eqref{eq:inv21}, \eqref{e:1} and \eqref{eq:Zes} together, we obtain for all $s \in (s_1,+\infty)$:
\begin{equation}
\label{eq:inv22}
\norm{e^{s\al(\cdot,0)}v(\cdot,0)}_{L^2(\Om)}^2 
\le C s^{-\f{1}{2}} \left( \norm{\pa_\nu v}_{L^2(\Sigma_0)}^2 + \norm{e^{s\al(\cdot,0)}\rho}_{L^2(\Om)}^2 + \norm{e^{s\al(\cdot,0)}A}_{L^2(\Om)}^2 + \norm{e^{s\al(\cdot,0)} \na\!\cdot\! A}_{L^2(\Om)}^2 \right).
\end{equation}
The last part of the proof is to lower estimate $\norm{e^{s\al(\cdot,0)}v(\cdot,0)}_{L^2(\Om)}$ in terms of $\norm{e^{s\al(\cdot,0)}\rho}_{L^2(\Om)}$, $\norm{e^{s\al(\cdot,0)}A}_{L^2(\Om)}$ and $\norm{e^{s\al(\cdot,0)} \na\!\cdot\! A}_{L^2(\Om)}$. To do that, we refer once more to the third line of \eqref{sys:maglinv}, giving
\begin{equation}
\label{eq:inv23}
e^{s\al(\cdot,0)}v(\cdot,0) =
-e^{s\al(\cdot,0)} \left(2 A\!\cdot\!\na u_0 + i \left( \rho + (A_1+A_2)\!\cdot\! A - i \na\!\cdot\! A \right) u_0\right),\
s \in \R,
\end{equation}
and we proceed by choosing $n+1$ times the initial state $u_0$ suitably, as described below.

\noindent{\it First choice}. We set $u_0=u_0^0$ where the function $u_0^0 : \Omega \to \Bbb{C}$ is constant and non-zero. Thus there exists $r_0 \in (0,+\infty)$ such that $\abs{u_0^0(x)} = r_0$ for all $x \in \Om$, and \eqref{eq:inv22}-\eqref{eq:inv23} then yield
\begin{equation}
\begin{aligned}
\label{eq:inv24}
& \norm{e^{s\al(\cdot,0)}(\rho + (A_1+A_2)\!\cdot\! A - i \na\!\cdot\! A)}_{L^2(\Om)}^2\\
\le & C\left(\norm{\pa_\nu v^0}_{L^2(\Sigma_0)}^2+ s^{-\f{1}{2}} \left( \norm{e^{s\al(x,0)}\rho}_{L^2(\Om)}^2 + \norm{e^{s\al(x,0)}A}_{L^2(\Om)}^2 + \norm{e^{s\al(\cdot,0)} \na\!\cdot\! A}_{L^2(\Om)}^2 \right) \right),\ s \in (s_1,+\infty),
\end{aligned}
\end{equation}
where $v^0:=\pa_t(u_1^0-u_2^0)$ and each $u_k^0$, $k=1,2$, is the solution to \eqref{sys:ori} associated with $u_0=u_0^0$, $\rho_0=\rho_k$ and $A_0=A_k$. \\

\noindent{\it Second choice.} We choose $n$ functions $u_0^j : \Om \to \R$, $j=1,\ldots,n$, such that the matrix $U_0^* U_0$, where $U_0:=\left( \pa_i u_0^j \right)_{1 \le i,j \le n}$ and $U_0^*$ denotes the Hermitian conjugate matrix to $U_0$, is strictly positive definite, i.e. such that
\begin{equation}
\label{e:4}
\exists r_0 \in (0,+\infty),\ \abs{U_0 \xi} \ge r_0 \abs{\xi},\ \xi \in \Bbb{C}^n.
\end{equation}
For each $j=1,\ldots,n$, we apply the well-known estimate 
\begin{equation}
\label{e:4b}
\abs{\xi+\zeta}^2 \ge \f{1}{2} \abs{\xi}^2 - \abs{\zeta}^2,\ \xi, \zeta \in \Bbb{C}^n,
\end{equation}
with $\xi=2 e^{s\al(\cdot,0)} A\!\cdot\!\na u_0^j$ and $\zeta= i e^{s\al(\cdot,0)} \left( \rho + (A_1+A_2)\!\cdot\! A - i \na\!\cdot\! A \right) u_0$ to \eqref{eq:inv23} and find
\begin{align*}
\norm{e^{s\al(\cdot,0)}v^j(\cdot,0)}_{L^2(\Om)}^2 \ge 2 \norm{e^{s\al(\cdot,0)} A\!\cdot\!\na u_0^j}_{L^2(\Om)}^2 - \norm{u_0^j}_{L^\infty(\Om)}^2 \norm{e^{s\al(\cdot,0)}(\rho + (A_1+A_2)\!\cdot\! A - i \na\!\cdot\! A)}_{L^2(\Om)}^2
\end{align*}
for all $s \in \R$. Here we have set $v^j:=\pa_t (u_1^j-u_2^j)$ where $u_k^j$, $k=1,2$, denotes the solution to \eqref{sys:ori} associated with $u_0=u_0^j$, $\rho=\rho_k$ and $A=A_k$. Summing up the above estimate over $j \in \{1,\ldots,n \}$ then yields 
\begin{align*}
\sum_{j=1}^n\norm{e^{s\al(\cdot,0)}v^j(\cdot,0)}_{L^2(\Om)}^2 &\ge 
2\norm{e^{s\al(\cdot,0)} U_0 A}_{L^2(\Om)}^2 - C\norm{e^{s\al(\cdot,0)}(\rho + (A_1+A_2)\!\cdot\! A - i\na\!\cdot\! A)}_{L^2(\Om)}^2\\
&\ge 2r_0^2\norm{e^{s\al(\cdot,0)}A}_{L^2(\Om)}^2 - C\norm{e^{s\al(\cdot,0)}(\rho + (A_1+A_2)\!\cdot\! A - i\na\!\cdot\! A)}_{L^2(\Om)}^2,\ s \in \R,
\end{align*}
in virtue of \eqref{e:4} and the identity $e^{s\al(\cdot,0)} U_0 A=U_0(e^{s\al(\cdot,0)} A)$. From this and \eqref{eq:inv22} it then follows for all $s \in (s_1,+\infty)$ that
\begin{align*}
\norm{e^{s\al(\cdot,0)}A}_{L^2(\Om)}^2 &\le C \left( \sum_{j=1}^n \norm{\pa_\nu v^j}_{L^2(\Sigma_0)}^2 + s^{-\f{1}{2}} \left( \norm{e^{s\al(\cdot,0)}\rho}_{L^2(\Om)}^2 + \norm{e^{s\al(\cdot,0)}A}_{L^2(\Om)}^2 + \norm{e^{s\al(\cdot,0)} \na\!\cdot\! A}_{L^2(\Om)}^2 \right) \right. \\
&\hspace{1cm} \left. +  \norm{e^{s\al(\cdot,0)}(\rho + (A_1+A_2)\!\cdot\! A - i \na\!\cdot\! A)}_{L^2(\Om)}^2 \right)
\end{align*}
and consequently 
\begin{equation}
\label{eq:inv25}
\norm{e^{s\al(\cdot,0)}A}_{L^2(\Om)}^2 \le C \left( \sum_{j=0}^n \norm{\pa_\nu v^j}_{L^2(\Sigma_0)}^2 + s^{-\f{1}{2}}\left(\norm{e^{s\al(\cdot,0)}\rho}_{L^2(\Om)}^2 + \norm{e^{s\al(\cdot,0)}A}_{L^2(\Om)}^2 + \norm{e^{s\al(\cdot,0)} \na\!\cdot\! A }_{L^2(\Om)}^2 \right) \right),
\end{equation}
from \eqref{eq:inv24}.

Further, by combining the following basic inequality
\begin{align*}
&\norm{e^{s\al(\cdot,0)}(\rho + (A_1+A_2)\!\cdot\! A -i \na\!\cdot\! A)}_{L^2(\Om)}^2 \\
 \ge & \f{1}{2}\norm{e^{s\al(\cdot,0)}(\rho -i \na\!\cdot\! A)}_{L^2(\Om)}^2 - \left( \norm{A_1}_{L^\infty(\Om)} + \norm{A_2}_{L^\infty(\Om)} \right)^2 \norm{e^{s\al(\cdot,0)}A}_{L^2(\Om)}^2,\ s \in \R,
\end{align*}
arising from \eqref{e:4b} where $\xi=e^{s \al(\cdot,0)}(\rho -i \na\!\cdot\! A)$ and $\zeta=e^{s\al(\cdot,0)}(A_1+A_2) \cdot A$,
with \eqref{eq:inv24} and \eqref{eq:inv25}, we get for all $s \in (s_1,+\infty)$ that
\begin{align}
& \norm{e^{s\al(\cdot,0)}(\rho -i \na\!\cdot\! A)}_{L^2(\Om)}^2 \nonumber \\
\le &  C \left( \sum_{j=0}^n \norm{\pa_\nu v^j}_{L^2(\Sigma_0)}^2 + s^{-\f{1}{2}} \left(\norm{e^{s\al(\cdot,0)}\rho}_{L^2(\Om)}^2 + \norm{e^{s\al(\cdot,0)}A}_{L^2(\Om)}^2 + \norm{e^{s\al(\cdot,0)} \na\!\cdot\! A}_{L^2(\Om)}^2 \right) \right).
\label{eq:inv26}
\end{align}
Having established \eqref{eq:inv26}, we turn now to estimating $\norm{e^{s\al(\cdot,0)}\rho}_{L^2(\Om)}$ and 
$\norm{e^{s\al(\cdot,0)} \na\!\cdot\! A}_{L^2(\Om)}$ with respect to the right hand side of \eqref{eq:inv25}. Let us first notice that we have
$\norm{e^{s\al(\cdot,0)}(\rho -i \na\!\cdot\! A)}_{L^2(\Om)}^2 = \norm{e^{s\al(\cdot,0)}\rho}_{L^2(\Om)}^2 + \norm{e^{s\al(\cdot,0)}(\na\!\cdot\! A)}_{L^2(\Om)}^2$ whenever the function $\rho$ is real-valued, in which case \eqref{eq:inv26} yields 
\begin{align}
& \norm{e^{s\al(\cdot,0)} \rho}_{L^2(\Om)}^2 +   \norm{e^{s\al(\cdot,0)} \na\!\cdot\! A}_{L^2(\Om)}^2 \nonumber \\
\le &  C \left( \sum_{j=0}^n \norm{\pa_\nu v^j}_{L^2(\Sigma_0)}^2 + s^{-\f{1}{2}} \left(\norm{e^{s\al(\cdot,0)}\rho}_{L^2(\Om)}^2 + \norm{e^{s\al(\cdot,0)}A}_{L^2(\Om)}^2 + \norm{e^{s\al(\cdot,0)} \na\!\cdot\! A}_{L^2(\Om)}^2 \right) \right),
\label{e:5}
\end{align}
for all $s \in (s_1,+\infty)$. In the general case where $\rho : \Om \to \Bbb{C}$, we combine the inequality $\abs{\na\!\cdot\! A} \le nM \abs{A}$ in $\Om$, arising from \eqref{h-em-b}, with \eqref{eq:inv25}, and get that
\begin{align*}
& \norm{e^{s \al(\cdot,0)} \na\!\cdot\! A}_{L^2(\Om)}^2 \\ 
\le & C \left( \sum_{j=0}^n \norm{\pa_\nu v^j}_{L^2(\Sigma_0)}^2 + s^{-\f{1}{2}} \left(\norm{e^{s\al(\cdot,0)}\rho}_{L^2(\Om)}^2 + \norm{e^{s\al(\cdot,0)}A}_{L^2(\Om)}^2 + \norm{e^{s\al(\cdot,0)} \na\!\cdot\! A}_{L^2(\Om)}^2 \right) \right), 
\end{align*}
for all $s \in (s_1,+\infty)$. This and \eqref{eq:inv26} yield \eqref{e:5} since $\norm{e^{s\al(\cdot,0)} \rho}_{L^2(\Om)} \le \norm{e^{s\al(\cdot,0)}(\rho -i \na\!\cdot\! A)}_{L^2(\Om)} + \norm{e^{s\al(\cdot,0)} \na\!\cdot\! A}_{L^2(\Om)}$
for all $s \in \R$. 

Now, putting \eqref{eq:inv25} and \eqref{e:5} together, we find that
\begin{align*}
& \norm{e^{s\al(\cdot,0)}\rho}_{L^2(\Om)}^2 + \norm{e^{s\al(\cdot,0)}A}_{L^2(\Om)}^2 + \norm{e^{s\al(\cdot,0)} \na\!\cdot\! A}_{L^2(\Om)}^2 \\ 
\le & C \left( \sum_{j=0}^n \norm{\pa_\nu v^j}_{L^2(\Sigma_0)}^2 + s^{-\f{1}{2}} \left(\norm{e^{s\al(\cdot,0)}\rho}_{L^2(\Om)}^2 + \norm{e^{s\al(\cdot,0)}A}_{L^2(\Om)}^2 + \norm{e^{s\al(\cdot,0)} \na\!\cdot\! A}_{L^2(\Om)}^2 \right) \right), 
\end{align*}
for all $s \in (s_1,+\infty)$. As a consequence there exists $s_2 \in [s_1,+\infty)$ so large that
$$
\norm{e^{s\al(\cdot,0)}\rho}_{L^2(\Om)}^2 + \norm{e^{s\al(\cdot,0)}A}_{L^2(\Om)}^2 + \norm{e^{s\al(\cdot,0)} \na\!\cdot\! A}_{L^2(\Om)}^2 \le C \sum_{j=0}^n \norm{\pa_\nu v^j}_{L^2(\Sigma_0)}^2,\ s \in (s_2,+\infty).
$$
Finally, this and \eqref{e:0} yield the stability estimate of Theorem  \ref{thm:stb2}.

\subsection{Proof of Theorem \ref{thm:stb3}}
\label{sec-stb3}
We stick with the notations of Subsection \ref{sec-stb2}.
The beginning of the proof of Theorem \ref{thm:stb3} follows the same path as the one of Theorem \ref{thm:stb2}, establishing the two estimates \eqref{eq:inv21} and \eqref{e:1}. They yield existence of a sufficiently large parameter $s_1 \in (0,+\infty)$ such that the following inequality
\begin{eqnarray}
& & \norm{e^{s\al(\cdot,0)} v(\cdot,0)}_{L^2(\Om)}^2  \nonumber \\
& \le &  Cs^{-\f{3}{2}} \left( \norm{\pa_\nu v}_{L^2(\Sigma_0)}^2 + \norm{e^{s\al(\cdot,0)}A}_{L^2(\Om)}^2 + \norm{e^{s \al(\cdot,0)} \na\!\cdot\! A }_{L^2(\Om)}^2 + s Z(v(\cdot,0)) \right), \label{s:0}
\end{eqnarray}
holds for all $s \in (s_1,+\infty)$. But, in the framework of Theorem \ref{thm:stb3} where none of the three assumptions \eqref{h-em-a}, \eqref{h-em-b} and \eqref{h-em-c} required by Theorem \ref{thm:stb2} is fulfilled, a more careful analysis is needed for majorizing the remaining term $Z(v(\cdot,0))$.

To this end we put $S:=A_1+A_2$ and we recall from Subsection \ref{sec-stb2} that $v(\cdot,0)=-2 A\!\cdot\!\na u_0 - i \left( \rho + S \cdot A - i \na\!\cdot\! A \right) u_0$ and that
$\nabla v(\cdot,0)= -(w_0 + z_0)$ with $w_0:= 2 \mathbb{D}^2_{u_0} A + (\nabla \cdot A) \nabla u_0 + i \left( \rho+ S \cdot A \right) \nabla u_0 + i \mathbb{J}_{S} A u_0 $ and $z_0:=2 \mathbb{J}_{A} \nabla u_0 + u_0 \nabla (\nabla \cdot A) + i \left( \nabla \rho+ \mathbb{J}_{A} S \right) u_0$. Thus,
$$
\ov{v(\cdot,0)} \nabla v(\cdot,0) -v(\cdot,0) \nabla \ov{v(\cdot,0)}=  2 i \left( \Im ( v(\cdot,0) \ov{w_0}) + \Im (v(\cdot,0) \ov{z_0}) \right) 
$$
so we may deduce from \eqref{defZ} that
\begin{align}
Z((v(\cdot,0)) & \le C \left( \int_{\Omega} e^{2s \al(x,0)} \abs{v(x,0) \ov{w_0(x)}} dx  + \int_{\Omega} e^{2s \al(x,0)} \abs{\Im (v(x,0) \ov{z_0(x)})}  dx \right) \nonumber \\
& \le C \left(  \norm{e^{s \al(x,0)} v(\cdot,0)}_{L^2(\Om)}^2 + \norm{e^{s \al(x,0)} w_0}_{L^2(\Om)}^2 + \int_{\Omega} e^{2s \al(x,0)} \abs{\Im (v(x,0) \ov{z_0(x)})}  dx\right),\ s \in \R. 
\label{s:1}
\end{align}
Since the two first terms in the right hand side of \eqref{s:1} are treated by \eqref{e:2} and \eqref{e:2b},  we may focus on the analysis of the third one. Actually, the imaginary part of $v(\cdot,0) \ov{z_0}$ decomposes into the sum of three terms, i.e. 
\begin{equation}
\label{s:1b}
\Im \left( v(\cdot,0) \ov{z_0} \right)  = \sum_{k=1}^3 \ell_k,
\end{equation}
with
$\ell_1 := \Im \left( \left( 2 A \cdot \nabla u_0 + (\nabla \cdot A) u_0 \right) \left( 2 \mathbb{J}_{A} \nabla \ov{u_0} + \ov{u_0} \nabla ( \nabla \cdot A) \right) \right)$,
$\ell_2 := \Im \left( (\rho + S \cdot A) (\nabla \ov{\rho} + \mathbb{J}_{A} S) \right) \abs{u_0}^2$ and
$\ell_3 :=\Re \left( (\rho + S \cdot A)u_0 (2 \mathbb{J}_{A} \nabla \ov{u_0} + \ov{u_0} \nabla ( \nabla \cdot A) )
- ( 2 A \cdot \nabla u_0 + (\nabla \cdot A) u_0 )  \ov{u_0} (\nabla \ov{\rho} + \mathbb{J}_{A} S)\right)$.
Since $A \in \R^n$, by assumption, we choose $u_0$ to be either real-valued or purely imaginary on $\Omega$, in such a way that we have
\begin{equation}
\label{s:2}
\ell_1(x)=0,\ x \in \Om.
\end{equation}
Further, as $\rho(x)=0$ and $S(x) \cdot A(x)= \abs{A_1(x)}^2 - \abs{A_2(x)}^2=0$ for a.e. $x \in \Omega$, it is apparent that
\begin{equation}
\label{s:3}
\ell_2(x) =0,\ x \in \Om,
\end{equation}
and $\ell_3=- ( 2 A \cdot \nabla u_0 + (\nabla \cdot A) u_0 )  \ov{u_0}  \mathbb{J}_{A} S$. Moreover, for all $i=1,\ldots,n$ and a.e. $x \in \Omega$, it holds true that
$$ 0  = \partial_i  \left( \abs{A_1(x)}^2 - \abs{A_2(x)}^2 \right) =  \partial_i \sum_{j=1}^n s_j(x) a_j(x) = \sum_{j=1}^n \left( \partial_i s_j(x) \right) a_j(x) + \sum_{j=1}^n \left( \partial_i a_j(x) \right) s_j(x),
$$
where we used the notations $S=(s_j)_{1 \le j \le n}$ and $A=(a_j)_{1 \le j \le n}$. Therefore, we have
$\mathbb{J}_{S} A + \mathbb{J}_{A} S=0$ and consequently $\ell_3=  (2 A \cdot \nabla u_0 + (\nabla \cdot A) u_0 )  \ov{u_0}  \mathbb{J}_{S} A$. This, \eqref{s:1b}, \eqref{s:2} and \eqref{s:3} yield
$$
\int_{\Om} e^{2s \al(x,0)} \abs{\Im ( v(x,0) \ov{z_0}(x) )} dx \leq C \left( \norm{e^{s \al(\cdot,0)} A}^2 + \norm{e^{s \al(\cdot,0)} \nabla \cdot A}^2 \right),\ s \in \R,
$$
hence
$$
Z((v(\cdot,0)) \le C \left( \norm{e^{s \al(\cdot,0)} A}^2 + \norm{e^{s \al(\cdot,0)} \nabla \cdot A}^2 \right),\
s \in \R,
$$
by \eqref{e:2}-\eqref{e:2b} and \eqref{s:1}. Now, putting this together with \eqref{s:0} we obtain for all $s \in (s_1,+\infty)$ that
\begin{equation}
 \label{s:3}
 \norm{e^{s\al(\cdot,0)} v(\cdot,0)}_{L^2(\Om)}^2  \le C s^{-\frac{1}{2}} \left( \norm{\pa_\nu v}_{L^2(\Sigma_0)}^2 + \norm{e^{s\al(\cdot,0)}A}_{L^2(\Om)}^2 + \norm{e^{s \al(\cdot,0)} \na\!\cdot\! A }_{L^2(\Om)}^2 \right).
\end{equation}
The rest of the proof follows the same lines as in the derivation of \eqref{eq:inv24} and \eqref{eq:inv25}. Namely, we first choose $u_0=u_0^0$ in \eqref{eq:inv23}, where $u_0^0(x)=r_0$ for some $r_0 \in \R \setminus \{ 0 \}$ and a.e. $x \in \Om$. We get that $e^{s\al(\cdot,0)} v(\cdot,0)=i e^{s\al(\cdot,0)} r_0 \nabla \cdot A$ for every $s \in \R$, hence
\begin{equation}
\label{s:4}
\norm{e^{s\al(\cdot,0)} \na\!\cdot\! A}_{L^2(\Om)}^2
\le C s^{-\frac{1}{2}} \left( \norm{\pa_\nu v^0}_{L^2(\Sigma_0)}^2+ \norm{e^{s\al(x,0)}A}_{L^2(\Om)}^2 + \norm{e^{s\al(\cdot,0)} \na\!\cdot\! A}_{L^2(\Om)}^2 \right),\ s \in (s_1,+\infty),
\end{equation}
according to \eqref{s:3}. Here $v^0:=\partial_t(u_1^0-u_2^0)$, where $u_j^0$ denotes the solution to \eqref{sys:ori} with $\rho_0=0$ and $A_0=A_j$ for $j=1,2$. Similarly, we choose $n$ real-valued functions $u_0^j$, $j=1,\ldots,n$, fulfilling \eqref{e:4}, take $u_0=u_0^j$ in \eqref{eq:inv23}, and argue as in the derivation of \eqref{eq:inv25} by substituting \eqref{s:4} for \eqref{eq:inv24}: we obtain that
$$
\norm{e^{s\al(\cdot,0)}A}_{L^2(\Om)}^2 \le C \left( \sum_{j=0}^n \norm{\pa_\nu v^j}_{L^2(\Sigma_0)}^2 + s^{-\f{1}{2}}\left(\norm{e^{s\al(\cdot,0)}A}_{L^2(\Om)}^2 + \norm{e^{s\al(\cdot,0)} \na\!\cdot\! A}_{L^2(\Om)}^2 \right) \right),\ s \in (s_1,+\infty),
$$
where $v^j:=\partial_t(u_1^j-u_2^j)$ and $u_k^j$, $k=1,2$, is the solution to \eqref{sys:ori} with $(\rho_0,A_0)=(0,A_k)$.
This and \eqref{s:4} yield the desired result.


\end{document}